\documentclass[a4paper]{amsart}
\usepackage{amsmath}
\usepackage{amsthm}
\usepackage{amssymb}
\usepackage{graphicx}
\usepackage{amsfonts}
\usepackage[applemac]{inputenc}
\usepackage{mathrsfs}
\usepackage{multirow}
\usepackage{pdfsync}
\usepackage{color}
\usepackage{bm}
\usepackage{mathdots}

\theoremstyle{plain}
\newtheorem{theorem}{Theorem}[section]

\theoremstyle{definition}

\newtheorem{remark}[theorem]{Remark}
\newtheorem*{remark*}{Remark}

\begin{document}
\title[]
{Multidomain spectral method for the Gauss hypergeometric function}

\author{S.~Crespo, M.~Fasondini, C.~Klein, N.~Stoilov,
C. Vall\'ee}

\address{School of Mathematics, Statistics and Actuarial Science, 
Sibson Building, Parkwood Road, 
University of Kent, 
Canterbury, Kent, CT2 7FS, 
United Kingdom.}
\email{M.Fasondini@kent.ac.uk}

\address{Institut de Math\'ematiques de Bourgogne
9 avenue Alain Savary, BP 47870, Universit\'e de 
Bourgogne-Franche-Comt\'e, 21078 Dijon Cedex, France}
\email{Christian.Klein@u-bourgogne.fr}

\address{Institut de Math\'ematiques de Bourgogne
9 avenue Alain Savary, BP 47870, Universit\'e de 
Bourgogne-Franche-Comt\'e, 21078 Dijon Cedex, France}
\email{Nikola.Stoilov@u-bourgogne.fr}

\date{\today}
\maketitle

\begin{abstract}
We present a multidomain spectral approach for Fuchsian ordinary differential 
equations in the particular case of the hypergeometric 
equation. Our hybrid approach uses Frobenius' method and Moebius 
transformations in the vicinity of each of the singular points of the 
hypergeometric equation, which leads to a natural decomposition of 
the real axis into domains. In each domain, solutions to the 
hypergeometric equation are 
constructed via the well-conditioned ultraspherical spectral method. The solutions are 
matched at the domain boundaries to lead to a solution which is 
analytic on the whole compactified real line $\mathbb{R}\cup 
{\infty}$, except for the singular points and cuts of the Riemann 
surface on which the solution is defined. The solution is further extended to the whole 
Riemann sphere by using the same approach for ellipses enclosing the 
singularities. The hypergeometric equation is solved on the ellipses 
with the boundary data from the real axis. This solution is continued 
as a harmonic function to the interior of the disk by solving the Laplace equation in polar 
coordinates with an optimal complexity Fourier--ultraspherical spectral method. 
\end{abstract}

\section{Introduction}\label{sec:intro}
Gauss' hypergeometric function $F(a,b,c,z)$ is arguably one of the most important 
classical transcendental functions in applications, see for instance 
\cite{Sea} and references therein. Entire  
 chapters are dedicated to it in various handbooks of mathematical functions such as the classical reference \cite{AS} and its modern 
reincarnation \cite{NIST}. It contains large classes of simpler transcendental 
functions as degeneracies, for example the Bessel functions. Despite its omnipresence 
in applications, numerical computation for a wide range of the 
parameters $a,b,c,z$ is challenging, see \cite{POP} for a 
comprehensive recent review with many references and a comparison of 
methods, and \cite{AKKW} for additional approaches to singular ODEs.    
This paper is concerned with the numerical evaluation of the 
hypergeometric function $F(a,b,c,z)$, treated  as a solution to a  
Fuchsian equation. This class of equations further includes examples such as the Lam\'e equation, 
see \cite{AS}, to which the method can be directly extended. The focus here is on the 
efficient computation of the hypergeometric function on the 
compactified real line $\mathbb{R}\cup\{\infty\}$ and on the Riemann 
sphere $\bar{\mathbb{C}}$, not just for individual values of $z$. The 
paper is intended as a proof of concept to study efficiently global 
(in the complex plane) 
solutions to singular ODEs (in these approaches, infinity is just a 
grid point and large values of the argument are treated as values in 
the vicinity of the origin) such as the Heun equation and Painlev\'e  
equations. 

The hypergeometric function can be 
defined in many ways, see for instance \cite{AS}. In this paper we 
construct it as the solution of the 
hypergeometric differential equation
\begin{equation}
    x(1-x)y''+(c-(1+a+b)x)y'-aby=0
    \label{hypergeom}
\end{equation}
with $F(a,b,c,0)=1$; here $a,b,c\in \mathbb{C}$ are constant with 
respect to $x\in \mathbb{R}$. 
Equation (\ref{hypergeom}) has regular (Fuchsian) singularities at $0$, $1$ and infinity. 
The Riemann symbol of the equation is given by
\begin{equation}
    P
    \begin{Bmatrix}
        0  & 1 & \infty & \\
        0 & 0 & a & z\\
        1-c & c-a-b & b&
    \end{Bmatrix}
    \label{symbol},
\end{equation}
where the three singularities are given in the first line. 
The second and third lines of the symbol (\ref{symbol}) 
give the exponents in the generalized series solutions of 
equation (\ref{hypergeom}): If $\xi$ is a local parameter near any of these 
singularities, following 
Frobenius' method (see again \cite{AS}), the general solution can be written  for sufficiently 
small $|\xi|$ in the form of generalized power series
\begin{equation}
    y = 
    \xi^{\kappa_{1}}\sum_{n=0}^{\infty}\alpha_{n}\xi^{n}+\xi^{\kappa_{2}}
    \sum_{m=0}^{\infty}\beta_{m}\xi^{m}
    \label{frobenius}
\end{equation}
if the difference between the constants $\kappa_{1}$ and $\kappa_{2}$ 
(corresponding to the second and third lines of the Riemann symbol 
(\ref{symbol}) respectively) is not integer; the constants (with 
respect to $\xi$) $\alpha_{n}$ 
and $\beta_{m}$ in (\ref{frobenius}) are given for $n>0$ and $m>0$ in 
terms of $\alpha_{0}$ and $\beta_{0}$, the last two being the only free constants 
in (\ref{frobenius}). It is known that generalized series of the form (\ref{frobenius}) have a radius of convergence equal to the 
minimal distance from the considered singularity to one of the other 
singularities of (\ref{hypergeom}). 

Note that logarithms may appear 
in the solution if the difference of the constants $\kappa_{1}$ and 
$\kappa_{2}$ is integer. We do not consider such cases here and concentrate on the generic case
\begin{equation}
    c, c-a-b, a-b\notin\mathbb{Z}.
    \label{generic}
\end{equation}
It is also well known that Moebius transformations (in other words elements of  $PSL(2, \mathbb{C})$)
\begin{equation}
    x\mapsto \frac{\alpha x+\beta}{\gamma x+\delta},\quad 
    \alpha,\beta,\gamma,\delta\in\mathbb{C}, \quad 
    \alpha\delta-\beta\gamma\neq 0,
    \label{moebius}
\end{equation}
transform one Fuchsian equation into another Fuchsian equation, the 
hypergeometric one (\ref{hypergeom}) into Riemann's differential 
equation, see \cite{AS}.   Moebius transformations can thus be used to 
map any singularity of a Fuchsian equation to 0.

The goal of this paper is to construct numerically the hypergeometric function $F(a,b,c,x)$
for generic arbitrary values of $a$, $b$, $c$, subject to
condition (\ref{generic}), and for arbitrary $x\in \bar{\mathbb{C}}$, and this 
with an efficient numerical approach showing spectral accuracy, i.e., 
an exponential decrease of the numerical error with the number of 
degrees of freedom used. We employ a hybrid strategy, that is we use Moebius transformations (\ref{moebius}) 
to map the considered singularity to 0, and we apply a change of the 
dependent variable so that the transformed solution is just 
the hypergeometric function in the vicinity of the origin, but with transformed 
values of the constants $a,b,c$. This is similar to 
Kummer's approach, see \cite{AS}, to express the solutions to the 
hypergeometric equation in different domains of the complex plane via the hypergeometric 
function near the origin. We thus obtain 3 
domains covering the complex plane, each of them centered at one of the three 
singular points of (\ref{hypergeom}). The dependent 
variable of the transformed equation (\ref{hypergeom}) is then 
transformed as $y\mapsto \xi^{\kappa_{i}}$, $i=1,2$, where the 
$\kappa_{i}$ are the exponents in (\ref{frobenius}). This 
transformation implies that we get yet another form of the Fuchsian 
equation which has a solution in terms of a power series. This 
solution will then be constructed numerically. 

This means that we 
solve one equation in the vicinity of $x=0$, and two each in each of 
the vicinities of $1$ and infinity. Instead of one form of the 
equation, we solve five $PSL(2,\mathbb{C})$  equivalent forms of 
(\ref{hypergeom}). 
This will be first done for real 
values of $x$. Since power series are in general 
slowly converging because of \emph{cancellation errors}, see for 
instance the discussions in \cite{POP},  we solve 
instead each of the 5 equivalent formulations of (\ref{hypergeom}) subject to the condition $y(\xi=0)=1$ (in an abuse of notation, we use the same 
symbol for the local variable $\xi$ and the dependent variable $y$ 
in all cases)   with spectral methods. Spectral methods are numerical methods for the global solution of differential equations that converge exponentially fast to analytic solutions. We shall use the efficient ultraspherical spectral method~\cite{US} which, as we shall see, can achieve higher accuracy than traditional spectral methods such as collocation methods because it is better conditioned. Solutions in each of the three domains are then matched (after 
multiplication with the corresponding factor $\xi^{\kappa}_{i}$, 
$i=1,2$) at the 
domain boundaries to the hypergeometric function constructed near 
$x=0$ to obtain a function which is $C^{1}$ at these boundaries (being a solution of the  
hypergeometric equation then guarantees the function is 
analytical if it is $C^{1}$ at the domain boundaries). Thus 
we obtain an analytic continuation of the hypergeometric function to 
the whole real line including infinity. 

Solutions to Fuchsian equations can be analytically continued as 
meromorphic functions to the whole complex plane (more precisely, 
they are meromorphic functions on a Riemann surface as detailed 
below). Since the Frobenius approach (\ref{frobenius}) is also 
possible with complex $x$, techniques similar to the approach for the 
real axis can be retained: we consider again three domains, each of them containing 
exactly one of the three singularities and the union of which is
covering the entire Riemann sphere $\bar{\mathbb{C}}$. On the boundary 
of each of these domains we solve the 5 equivalent forms of 
(\ref{hypergeom}) as on the real axis with boundary data obtained on 
$\mathbb{R}$. Then the holomorphic function corresponding to the 
solution in the interior of the studied domain is obtained by solving 
the Laplace equation with the data obtained on the boundary. 
The advantage of the Laplace equation is that it is not 
singular in contrast to the hypergeometric equation. 
It is solved by introducing polar coordinates $r,\phi$ in each of the 
domains and then using
 the ultraspherical spectral method in $r$ and a 
Fourier spectral approach in $\phi$ (the solution is periodic in 
$\phi$). Since the matching has already been done on the real axis, 
one immediately obtains the hypergeometric function on 
$\bar{\mathbb{C}}$ in this way. 

The paper is organized as follows: in section 2 we construct the 
hypergeometric function on the real axis. In section 3 the 
hypergeometric function is analytically continued to a meromorphic 
function on the Riemann sphere. In section 4 we consider examples for 
interesting values of $a,b,c,x$ in (\ref{hypergeom}). In section 5 we add some concluding 
remarks. 

\section{Numerical construction of the hypergeometric function on the 
real line}\label{sec2}

In this section we construct numerically the hypergeometric function 
on the whole compactified real line. To this end we introduce the 
following three domains:\\
domain~I: local parameter $x$, $x\in [-1/2,1/2]$,\\
domain~II: local parameter $t = 1-x$, $t\in [-1/2,1/2]$,\\
domain~III:  local parameter $s = -1/(x-1/2)$, $s\in [-1,1]$.\\
In each of these domains we apply transformations to the dependent 
variable such that the solutions are analytic functions on their domains. Since we shall approximate the solution using Chebychev polynomials that are defined on the unit interval $[-1, 1]$, we map the above intervals $[x_{l}, x_r]$ to $[-1,  1]$   via $x_{l}(1-\ell)/2+x_{r}(1+\ell)/2$, $\ell \in [-1, 1]$. For these equations we 
look for the unique solutions with $y(0)=1$ since the solutions are all hypergeometric functions (which is defined to be $1$ at the origin) with transformed values of the parameters, as showed by Kummer~\cite{AS}. This 
means we are always studying equations of the form   
\begin{equation}
    a_2(\ell)y''+ a_1(\ell)y'+ a_0(\ell)y=0, \qquad \ell \in [-1, 1], \qquad y(0) = 1
    \label{form},
\end{equation}
where $a_2(\ell)$, $a_1(\ell)$ and $a_0(\ell)$  are polynomials and $a_2(0) = 0$.   At the domain boundaries, the solutions are matched by a $C^{1}$ condition 
on the hypergeometric function which is thus uniquely determined (the 
hypergeometric equation then implies that the solution is in fact 
analytical if it is $C^{1}$ at the domain boundaries). To 
illustrate this procedure we use the fact that many hypergeometric 
functions can be given in terms of elementary functions, see for 
instance \cite{AS}. Here we consider the example
\begin{equation}
    F(a,b,c,x) = (1-x)^{-a}, \qquad a=-1/3,\: b=c=1/2.
    \label{example}
\end{equation}
 The triple $a,b, c$ is thus generic as per condition 
(\ref{generic}). More general examples are discussed 
quantitatively in section 4. 

\subsection{The ultraspherical (US) spectral method}


The recently introduced ultraspherical (US) spectral method~\cite{US} overcomes some of the weaknesses of traditional spectral methods such as dense and ill-conditioned matrices. The key idea underlying the US method is to change the basis of the solution expansion upon differentiation to the ultraspherical polynomials, which leads to sparse and well-conditioned matrices, as we briefly illustrate below.

Since the solutions we compute are analytic, we can express $y$ as a Chebychev series~\cite{trefethen2}:
\begin{equation}
   y  = \sum_{j = 0}^{\infty} y_j T_j(\ell), \qquad T_{j}(\ell) = \cos[j \arccos(\ell)], \qquad \ell \in [-1, 1].  \label{Chebseries}
\end{equation} 
The differentiation operators in the US method are based on the following relations involving the ultraspherical (or Gegenbauer) orthogonal polynomials, $C^{(\lambda)}_j(\ell)$:
\begin{equation*}
\frac{\mathrm{d} T_j}{\mathrm{d}\ell} = \begin{cases}
j C_{j-1}^{(1)} & j \geq 1 \\
0 & j = 0
\end{cases}, \qquad
\frac{\mathrm{d} C_j^{(\lambda)}}{\mathrm{d}\ell} = \begin{cases}
2\lambda C_{j-1}^{(\lambda + 1)} & j \geq 1 \\
0 & j = 0
\end{cases}, \qquad \lambda \geq 1.
\end{equation*}
Hence, differentiations of 
(\ref{Chebseries}) give
\begin{equation}
y'  = \sum_{j = 0}^{\infty}(j+1) y_{j+1} C^{(1)}_j(\ell) \quad \text{and} \quad y'' = 2\sum_{j = 0}^{\infty}(j+2) y_{j+2} C^{(2)}_j(\ell). \label{USders}
\end{equation}
Therefore, if we let $\bm{y}$ denote the (infinite) vector of Chebychev coefficients of $y$, then the coefficients of  the $\lbrace C^{(1)}_j \rbrace $ and $\lbrace C^{(2)}_j \rbrace $ expansions of $y'$ and $y''$ are given by $\mathcal{D}_1\bm{y}$ and $\mathcal{D}_2\bm{y}$, respectively, where 
\begin{equation*}
\mathcal{D}_1 = \left(
{\renewcommand{\arraystretch}{1.0}
\begin{array}{c c c c c}
0 & 1 &   &        &   \\
  &   & 2 &        &   \\
  &   &   & 3 &  \\
  &   &   &        & \ddots
\end{array}
}
  \right) \quad \text{and} \quad
\mathcal{D}_2 = 2\left(
{\renewcommand{\arraystretch}{1.0}
\begin{array}{c c c c c c}
0 & 0 & 2  &         &        &   \\
  &   &    &  3      &        &    \\
  &   &    &         & 4 &    \\
  &   &    &         &        & \ddots
\end{array}
}
  \right). 
\end{equation*}


Notice that the expansions in (\ref{Chebseries}) and (\ref{USders}) are expressed in different polynomial bases ($\lbrace T_j \rbrace $, $\lbrace C^{(1)}_j \rbrace $, $\lbrace C^{(2)}_j \rbrace $). The next steps in the US method are (i) substitute (\ref{Chebseries}) and (\ref{USders}) into the differential equation (\ref{form}) and perform the multiplications $a_2(\ell) y'' $,  $a_1(\ell) y' $ and $a_0(\ell) y$ in the $\lbrace C^{(2)}_j \rbrace $, $\lbrace C^{(1)}_j \rbrace $ and $\lbrace T_j \rbrace $  bases, respectively, and (ii) convert the expansions in the $\lbrace T_j \rbrace $ and $\lbrace C^{(1)}_j \rbrace $ bases to the $\lbrace C^{(2)}_j \rbrace $ basis.


Concerning step (i), consider the term $a_0(\ell)y$ and suppose that $a_0(\ell) = \sum_{j = 0}^{\infty}a_j T_j(\ell)$, then
\begin{equation*}
a_0(\ell)y =   \sum_{j = 0}^{\infty} a_j T_j(\ell) \sum_{j = 0}^{\infty} y_j T_j(\ell) = \sum_{j= 0}^{\infty}c_j T_j(\ell),
\end{equation*}
where~\cite{US}
\begin{equation}
c_j = \begin{cases}
a_0 y_0 + \frac{1}{2} \sum_{k = 1}^{\infty} a_k y_k & j = 0 \\
\frac{1}{2} \sum_{k = 0}^{j-1} a_{j-k} y_k + a_0 y_j + \frac{1}{2} \sum_{k = 1}^{\infty} a_{k} y_{k+j} + \frac{1}{2} \sum_{k = 0}^{\infty} a_{k+j} y_{k} & j \geq 1.  \label{multcoeffs}
\end{cases}
\end{equation}
Expressed as a multiplication operator on the Chebychev coefficients $\bm{y}$, (\ref{multcoeffs}) becomes  $\bm{c} = \mathcal{M}_0[a_0(\ell)] \bm{y}$, where $\mathcal{M}_0[a_0(\ell)]$ is a Toeplitz plus an almost Hankel operator given by
\begin{equation}
\mathcal{M}_0[a_0(\ell)] = \frac{1}{2}\left[
\left(  
\begin{array}{c c c c c}
2a_0 & a_1 & a_2 & a_3 & \cdots \\
a_1  &2a_0 & a_1 & a_2 & \ddots \\
a_2 & a_1  &2a_0 & a_1 & \ddots \\
a_3 &a_2   & a_1 &2a_0 & \ddots \\
\vdots & \ddots & \ddots & \ddots & \ddots 
\end{array}
\right)
+ 
\left(  
\begin{array}{c c c c c}
0 & 0 & 0 & 0 & \cdots \\
a_1  & a_2 & a_3 & a_4 & \iddots \\
a_2  & a_3 & a_4 & a_5 & \iddots \\
a_3  & a_4 & a_5 & a_6 & \iddots \\
\vdots & \iddots & \iddots & \iddots & \iddots 
\end{array}
\right)
\right].  \label{multop}
\end{equation}
For all the equations considered in this section, $a_0(\ell)$ is either a zeroth or first degree polynomial and hence $a_j = 0$ for $j > 0$ or $j > 1$. In the next section, $a_0(\ell)$ will be analytic (and entire), in which case $a_0(\ell)$ can be uniformly approximated to any desired accuracy by using only $m$ coefficients $a_j, j = 0, \ldots, m-1$ for sufficiently large $m$ ($m = 36$ will be sufficient for machine precision accuracy in the next section). Therefore if $n > m$, then the multiplication operator (\ref{multop}) is banded with bandwidth $m-1$ (i.e., $m-1$ nonzero diagonals on either side of the main diagonal). 

In a similar vein, the multiplication of the series $a_2(\ell) y'' $ and $a_1(\ell) y' $  can be expressed in terms of the multiplication operators $\mathcal{M}_2[a_2(\ell)]$ and $\mathcal{M}_1[a_1(\ell)]$ operating on the coefficients of $ y'' $ (in the $\lbrace C^{(2)}_j \rbrace $ basis) and $ y' $ (in the $\lbrace C^{(1)}_j \rbrace $ basis), i.e., $\mathcal{M}_2[a_2(\ell)] \mathcal{D}_2 \bm{y}$ and $\mathcal{M}_1[a_1(\ell)] \mathcal{D}_1 \bm{y}$. If we represent or approximate $a_2(\ell)$ and $a_1(\ell)$ with $m$ Chebychev coefficients, then, as with $\mathcal{M}_0[a_0(\ell)]$, the operators $\mathcal{M}_2[a_2(\ell)]$ and $\mathcal{M}_1[a_1(\ell)]$ are banded with bandwidth $m-1$ if $n > m$. The entries of these operators are given explicitly in~\cite{US}.

For step (ii), converting all the series to the $\lbrace C^{(2)}_j \rbrace $ basis, we use
\begin{equation*}
T_j = \begin{cases}
\frac{1}{2}\left( C_{j}^{(1)} - C_{j-2}^{(1)} \right) & j \geq 2 \\
\frac{1}{2} C_{1}^{(1)} & j = 1 \\
 C_{0}^{(1)} & j = 0 
\end{cases}, \qquad
C_j^{(1)} = \begin{cases}
\frac{1}{1 + j}\left( C_{j}^{(2)} - C_{j-2}^{(2)} \right) & j \geq 2 \\
\frac{1}{2} C_{1}^{(2)} & j = 1 \\
 C_{0}^{(2)} & j = 0 
\end{cases},
\end{equation*}
from which the operators for converting the coefficients of a series from the basis $\lbrace T_j \rbrace $ to $\lbrace C_j^{(1)} \rbrace $ and from $\lbrace C_j^{(1)} \rbrace $ to $\lbrace C_j^{(2)} \rbrace $ follow, respectively:
\begin{equation*}
\mathcal{S}_0 = \left(
{\renewcommand{\arraystretch}{1.0}
\begin{array}{c c c c c c}
1 &              &-\frac{1}{2}   &              &  &   \\
  & \frac{1}{2}  &               & -\frac{1}{2} &  &   \\
  &              & \frac{1}{2}   &              & -\frac{1}{2} &  \\
  &              &               & \ddots       & & \ddots
\end{array}
}
  \right), \:  \:
\mathcal{S}_1 = \left(
{\renewcommand{\arraystretch}{1.0}
\begin{array}{c c c c c c}
1 &              &-\frac{1}{3}   &              &  &   \\
  & \frac{1}{2}  &               & -\frac{1}{4} &  &   \\
  &              & \frac{1}{3}   &              & -\frac{1}{5} &  \\
  &              &               & \ddots       & & \ddots
\end{array}
}
  \right). 
\end{equation*}
Thus, the linear operator
\begin{equation*}
\mathcal{L} := \mathcal{M}_2[a_2(\ell)] \mathcal{D}_2 + \mathcal{S}_1\mathcal{M}_1[a_1(\ell)] \mathcal{D}_1 + \mathcal{S}_1 \mathcal{S}_0 \mathcal{M}_0[a_0(\ell)],  
\end{equation*}
operating on the Chebychev coefficients of the solution, i.e., $\mathcal{L}\bm{y}$, gives the coefficients of the differential equation (\ref{form}) in the $\lbrace C_j^{(2)} \rbrace $ basis.

The solution (\ref{Chebseries}) is approximated by the first $n$ terms in its Chebychev expansion,
\begin{equation*}
y \approx \widetilde{y}_n := \sum_{j = 0}^{n-1} y_j T_j(\ell),
\end{equation*}
that satisfies the condition $\widetilde{y}_n(0) = 1$. To obtain an $n\times n$ linear system for these coefficients, the $\infty \times \infty$ operator $\mathcal{L}$ operator needs to be truncated  using the $n \times \infty$ projection operator given by
\begin{equation*}
\mathcal{P}_n = (I_n,\bm{0}),
\end{equation*} 
where $I_n$ is the $n\times n$ identity matrix. The  $n-1 \times n$ truncation of $\mathcal{L}$ is $\mathcal{P}_{n-1}\mathcal{L}\mathcal{P}_n^{\top}$, which is complemented with the condition $\widetilde{y}_n(0) = 1$. Then the $n \times n$ system to be solved is 
\begin{equation}
 \left( 
 {\renewcommand{\arraystretch}{1.5}
\begin{array}{c c c c}
T_0(0) & T_1(0) & \cdots & T_{n-1}(0) \\
 &\multicolumn{2}{c}{\mathcal{P}_{n-1}\mathcal{L}\mathcal{P}_n^{\top}}  & 
\end{array}}
\right)
\left(
\begin{array}{c}
y_0 \\
y_1 \\
\vdots \\
 y_{n-1}
\end{array}
\right)
= 
\left(
\begin{array}{c}
1 \\
0 \\
\vdots \\
0
\end{array}
\right),
\label{USsystem}
\end{equation}
where $T_j(0) = \cos(j\pi/2)$. We construct the matrix in (\ref{USsystem}) by using the functions provided in the Chebfun Matlab package~\cite{Chebfun}. Chebfun is also an ideal environment for stably and accurately performing operations on the approximation $\widetilde{y}_n$ (e.g., evaluation (with barycentric interpolation~\cite{barycentric}) and differentiation, which we require in the following section when matching the solutions at the domain boundaries).

\subsection{Domain I}

\begin{figure}[htb!]
\hspace*{-0.0 cm}
\mbox{\includegraphics[width=0.5\textwidth]{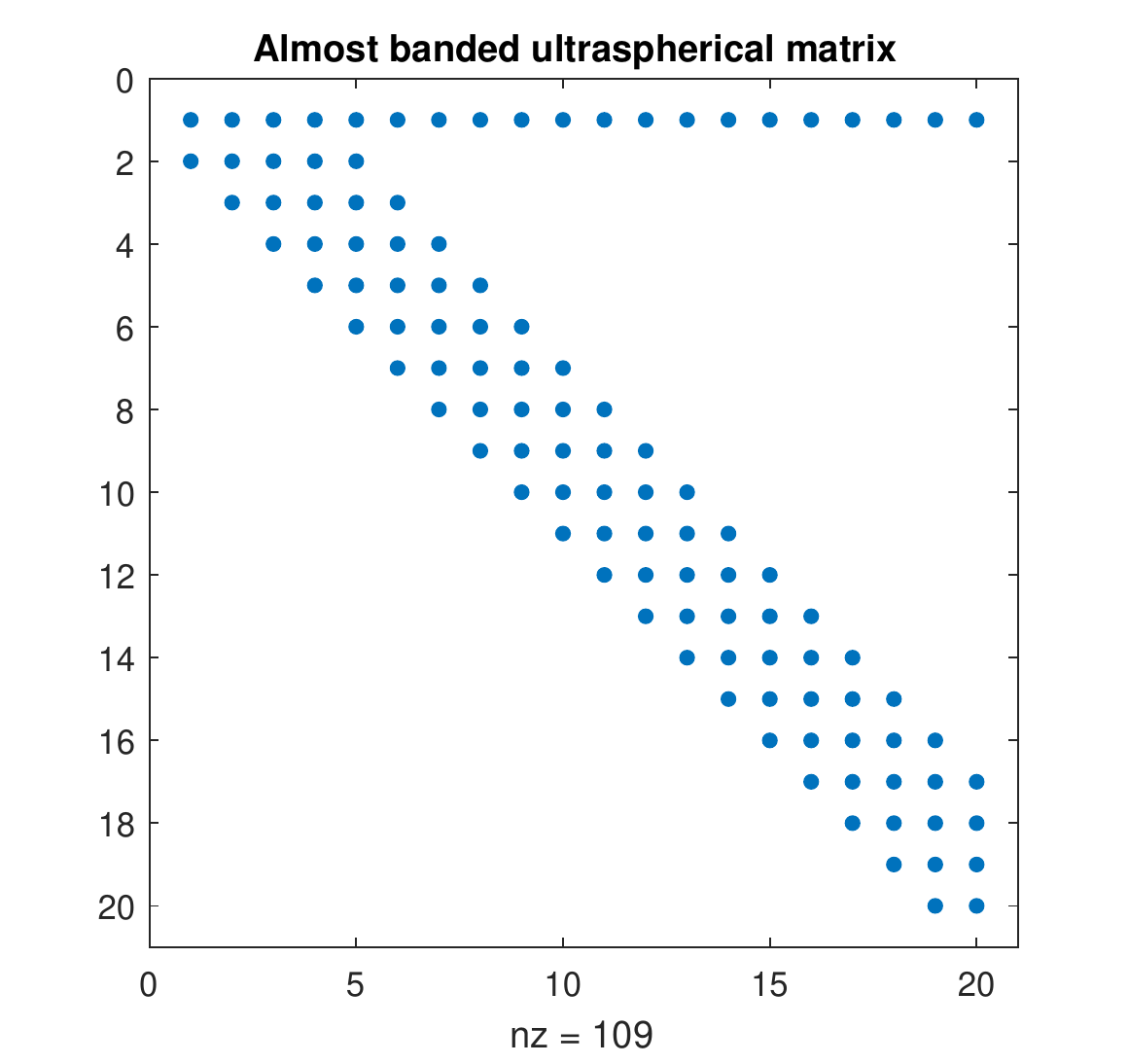}\includegraphics[width=0.5\textwidth]{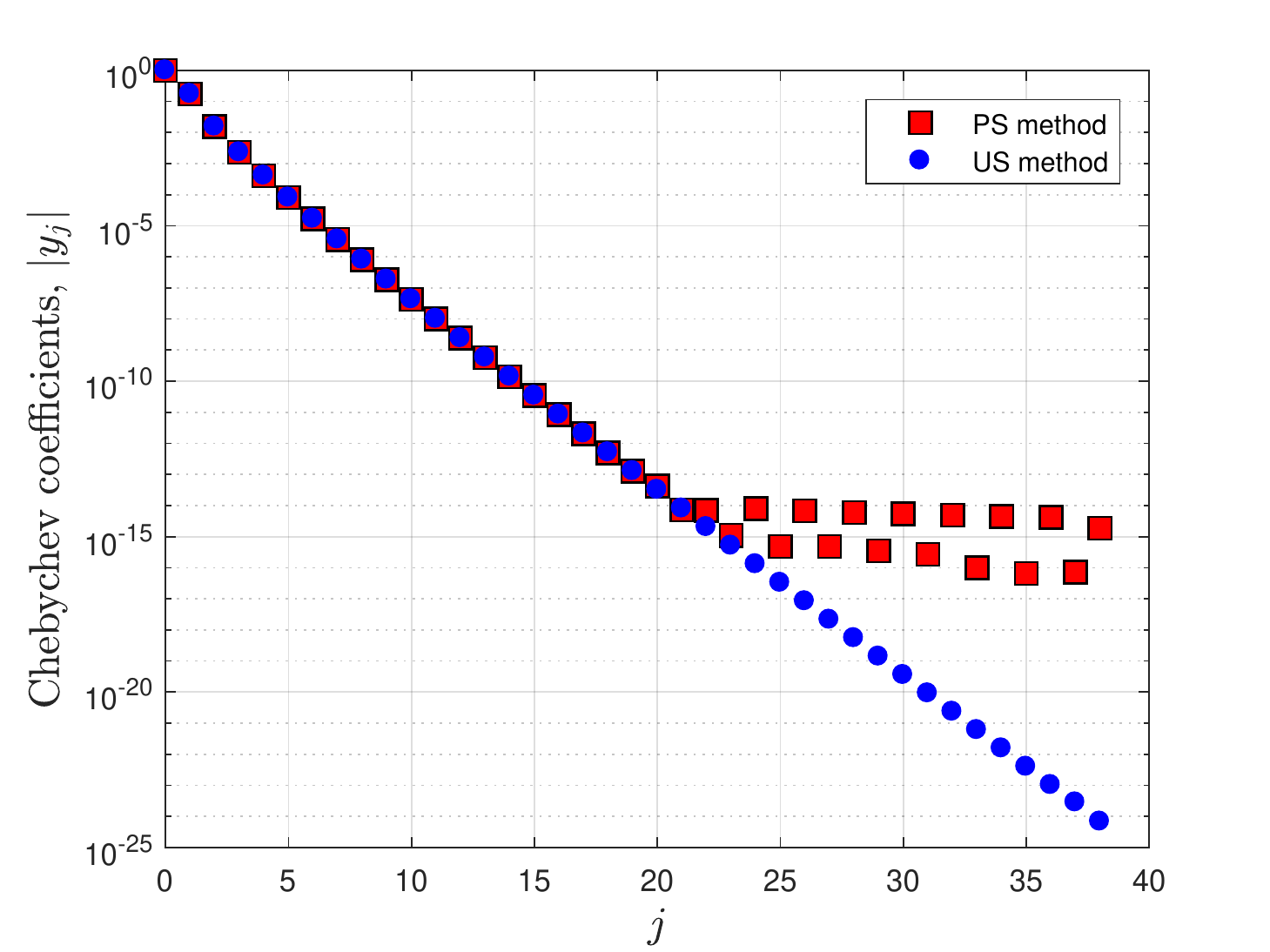}}
\hspace*{-0.0 cm}
\mbox{\includegraphics[width=0.5\textwidth]{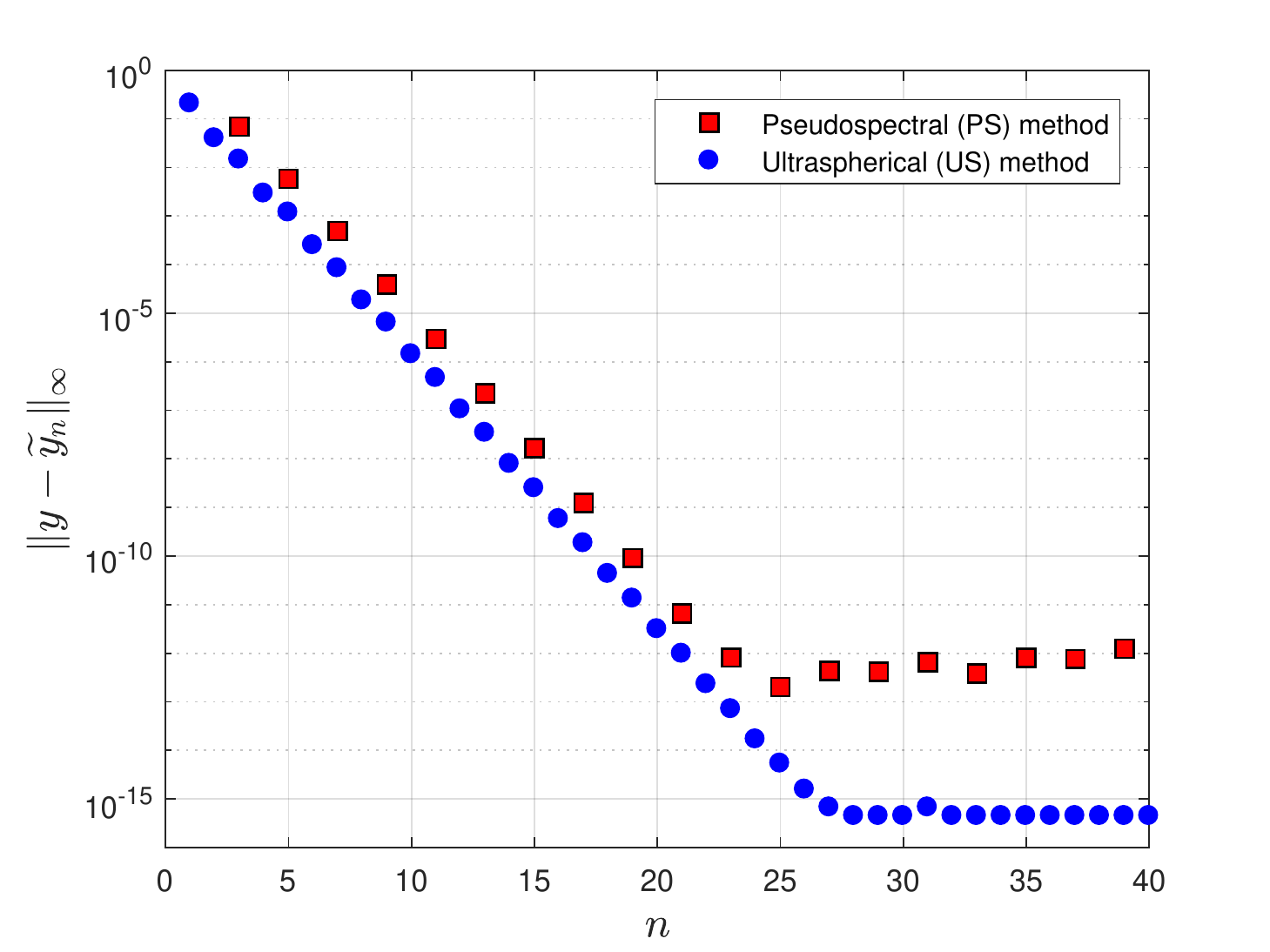}
\includegraphics[width=0.5\textwidth]{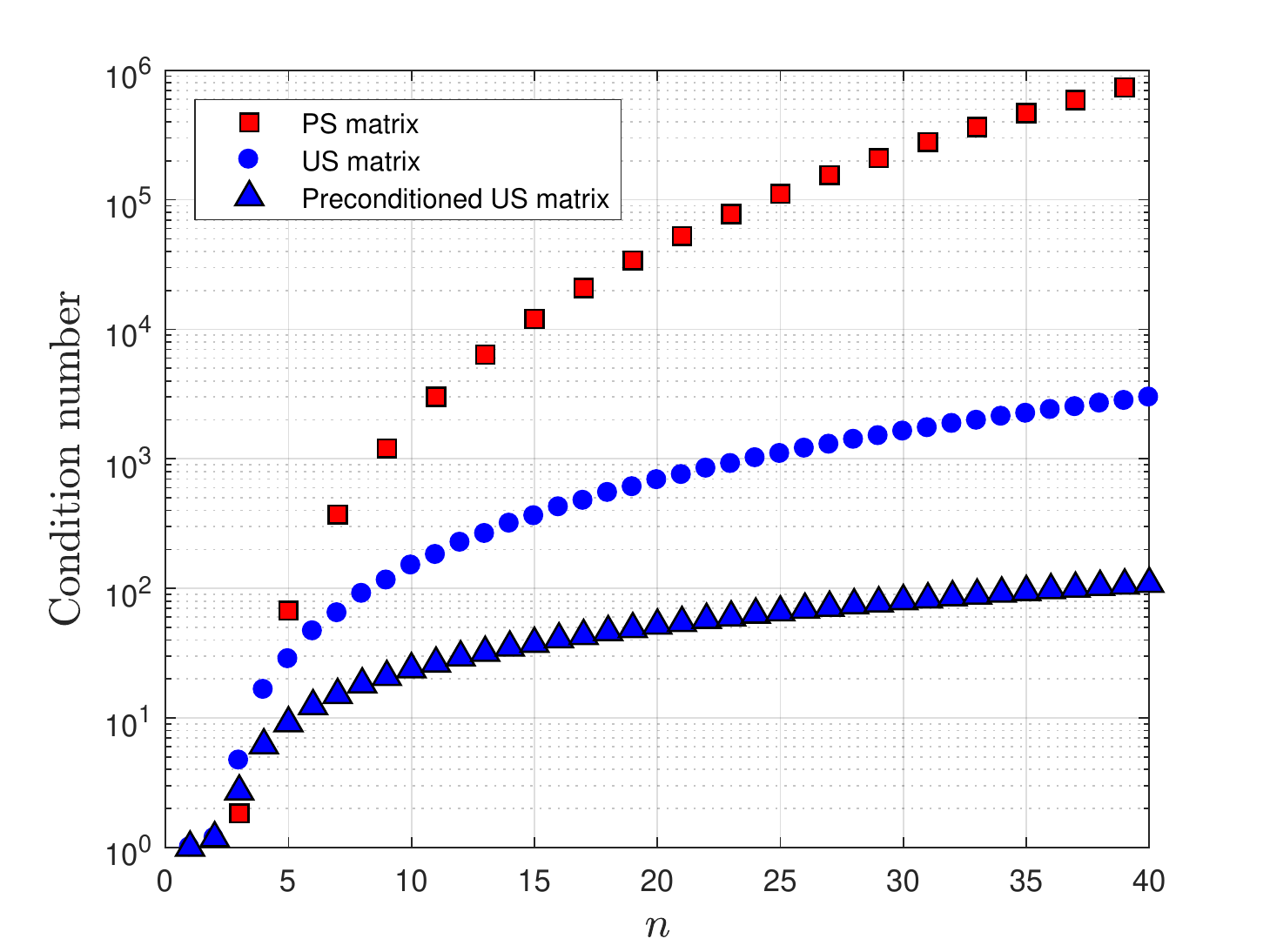}}
 \caption{The performance of the ultraspherical (US) spectral method and a Chebychev collocation, or pseudospectral (PS) method for the solution of the hypergeometric equation (\ref{hypergeom}) on the interval $[-1/2, 1/2]$ with $y(0) = 1$, $a = -1/3$ and $b = c = 1/2$. }
 \label{realI}
\end{figure}

Figure~\ref{realI} gives the results for solving (\ref{hypergeom}) for the test example (\ref{example}) on the interval $[-1/2, 1/2]$. 
 For comparison purposes, a Chebychev collocation, or pseudospectral (PS) method~\cite{trefethen} is also used. In contrast to the US spectral method in which the operators operate in coefficient space, in the PS method the matrices operate on the solution values at collocation points (e.g., the Chebychev points (of the second kind) $\cos(j\pi/n), j = 0, \ldots, n$). These matrices can be constructed in Matlab with Chebfun or the Differentiation Matrix Suite~\cite{WR}.

The top-left frame shows the almost banded structure of the matrix in the system (\ref{USsystem}), which can be solved in only $\mathcal{O}\left( m^2 n \right)$ operations using the adaptive QR method in~\cite{US}. Since this algorithm is not included in Matlab, we use the backslash command to solve the linear systems. The PS method, by comparison, yields dense linear systems. 

The top-right frame shows the magnitude of the $n=40$ Chebychev coefficients of the solution obtained with the US and PS methods. As expected, the magnitudes decrease exponentially with $n$ since the solution is analytic on the interval. 

The bottom-left frame shows the maximum error of the computed solution $\widetilde{y}_n$ on the interval $[-1/2, 1/2]$,  which can be accurately approximated in Chebfun. The solution converges exponentially fast to the exact solution for both methods, however, the US method stably achieves almost machine precision accuracy (on the order of $10^{-16}$) while the PS method reaches an accuracy of around $10^{-13}$ at $n=25$ and then the error increases slightly as $n$ is further increased.

The higher accuracy attainable by the US method and the numerical instability of the PS method are partly explained by the condition numbers of the matrices that arise in these methods (see the bottom-right frame). In~\cite{US} it is shown that, provided the equation has no singular points on the interval, the condition number of the US matrices grow linearly with $n$ and with preconditioning the condition number can be bounded for all $n$. By contrast, the condition numbers of collocation methods increase as $\mathcal{O}\left( n^{2N} \right)$, where $N$ is the order of the differential equation. Since equation (\ref{hypergeom}) has a singular point on the interval, we find different asymptotic growth rates of the condition numbers (by doing a least squares fit on the computed condition numbers): $\mathcal{O}\left(n\right)$, $\mathcal{O}\left(n^2\right)$
and $\mathcal{O}\left(n^{4.17}\right)$ for, respectively, preconditioned US matrices (using the diagonal preconditioner in~\cite{US}), the US matrices in (\ref{USsystem}) (with no preconditioner) and PS matrices.  We find that, as observed in~\cite{US}, the accuracy achieved by the US method is much better than the most pessimistic bound based on the condition number of the matrix---hardly any accuracy is lost despite condition numbers on the order of $10^3$. As pointed out in~\cite{US}, while a diagonal preconditioner decreases the condition number of the US matrix, solving a preconditioned system does not improve the accuracy of the solution if the system is solved using QR. This agrees with our experience that the accuracy obtained with Matlab's backslash solver does not improve if some digits are lost by the US method. Hence, all the numerical results reported in this paper were computed without preconditioning. 

%
%
%
%
%
%
%

\subsection{Domain II}

Next we address domain II with $x\in[0.5,1.5]$, where we use the 
local parameter $t=1-x$, in which (\ref{hypergeom}) reads 
\begin{equation}
    t(1-t)u''+(a+b+1-c-(1+a+b)t)u'-abu=0.
    \label{hypergeomt1}
\end{equation}
The solution corresponding 
to the exponent 0 in the symbol (\ref{symbol}) is constructed as in 
domain I with the US method.
There does not appear to be an 
elementary closed form of this solution for the studied example, which is plotted in Figure~\ref{real_sols_fig}.
 The results obtained for (\ref{hypergeomt1}), and also for the three equations remaining to be solved, (\ref{II}), (\ref{infv2}) and (\ref{infv3}), 
  are qualitatively the same as those in Figure~\ref{realI} and therefore we do not plot the results again.

%
%
%

The solution proportional to $t^{c-a-b}$ 
is constructed by writing $u=t^{c-a-b}\widetilde{u}(t)$. 
Equation 
(\ref{hypergeomt1}) implies for $\widetilde{u}$ the  equation 
\begin{equation}
    t(1-t)\widetilde{u}''+(c-a-b+1-t(2c-a-b+1))\widetilde{u}'-(b-c)(a-c)\widetilde{u}=0
    \label{II}.
\end{equation}
The solution to (\ref{II}) with $\widetilde{u}(0)=1$ is $\widetilde{u}=1$ since $b = c$ for the test problem and it is recovered exactly by the US method.

\begin{remark}\label{rem1}
    The appearance of the root 
$t^{c-a-b}$ in $u=t^{c-a-b} \widetilde{u}$ indicates that the solution as well as the hypergeometric 
function will in general not be single valued on $\mathbb{C}$, but on 
a Riemann surface. If the genericness condition (\ref{generic}) is 
satisfied, this surface will be compact. If this were not the case, 
logarithms could appear which are only single valued on a non-compact surface. 
To obtain a single valued function on a 
compact Riemann surface, monodromies have to be traced which can be 
numerically done as in \cite{RS}. This is beyond the scope of the 
present paper. Here we only construct the solutions to various 
equations which are entire and thus single valued on $\mathbb{C}$. 
The roots appearing in the representation of the hypergeometric 
function built from these single valued solutions 
are taken to be branched along the negative real axis. Therefore cuts 
may appear in the plots of the hypergeometric function.
\end{remark}

\subsection{Domain III}
For $x\sim \infty$, we use 
the local coordinate $s=-1/(x-1/2)$ with $s\in[-1,1]$. In this case we get 
for the hypergeometric equation
\begin{equation}
    \frac{s^{2}}{4}(s-2)(s+2)y''+sy'\left[\frac{s^{2}}{2}+\left(c-\frac{a+b+1}{2}\right)s
    +a+b-1\right]-aby=0
    \label{inf2}.
\end{equation}
Writing $y = s^{a}v$, we get for (\ref{inf2})
\begin{equation}
    \begin{split}
    &\frac{s}{4}(s-2)(s+2)v''+v'\left[(a+1)\frac{s^{2}}{2}+\left(c-\frac{a+b+1}{2}\right)s
    +b-a-1\right]\\
    &+a\left[c-\frac{a+b+1}{2}
    +\frac{s}{4}(a+1)\right]v=0.
    \end{split}
    \label{infv2}
\end{equation}

The hypergeometric equation (\ref{hypergeom}) is obviously invariant 
with respect to an interchange of $a$ and $b$. Thus we can 
always consider the case $\Re b>\Re a$. The solution of (\ref{inf2}) 
proportional to $s^{b}$  can be found by writing $y=s^{b}\widetilde{v}$ and exchanging $a$ and $b$ in 
(\ref{infv2}) ($b-a$ is not an integer because of (\ref{generic})), 
\begin{equation}
    \begin{split}
    &\frac{s}{4}(s-2)(s+2)\widetilde{v}''+\widetilde{v}'\left[(b+1)\frac{s^{2}}{2}+\left(c-\frac{a+b+1}{2}\right)s
    +a-b-1\right]\\
    &+b\left[c-\frac{a+b+1}{2}
    +\frac{s}{4}(b+1)\right]\widetilde{v}=0.
    \end{split}
    \label{infv3}
\end{equation}

\begin{figure}[htb!]
\includegraphics[width=0.525\textwidth]{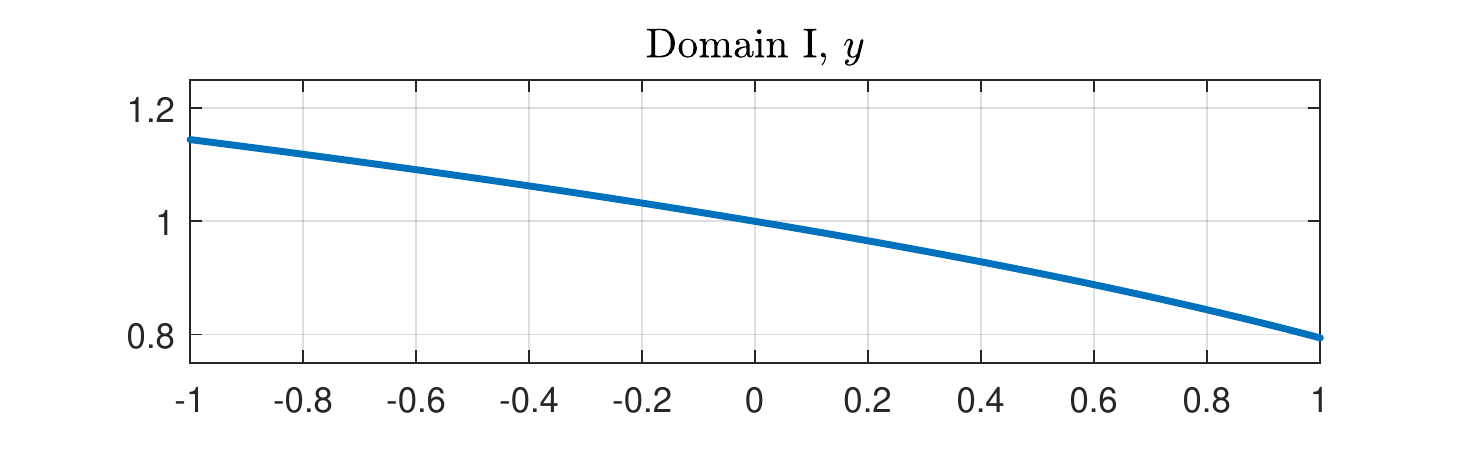}\\
\hspace*{-0.25 cm}
\mbox{\includegraphics[width=0.525\textwidth]{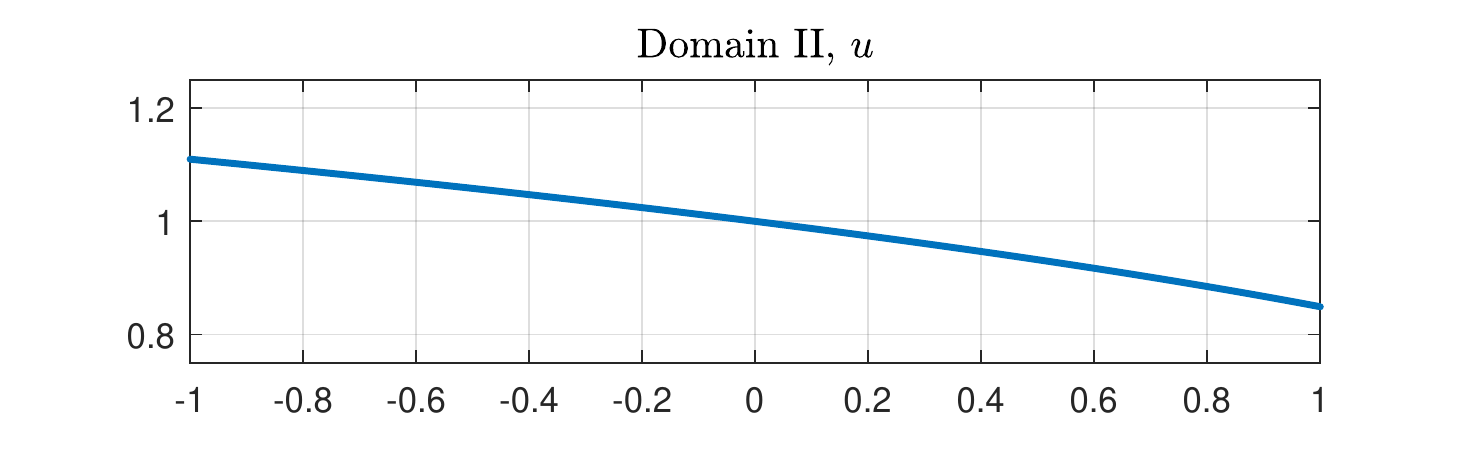}
\includegraphics[width=0.525\textwidth]{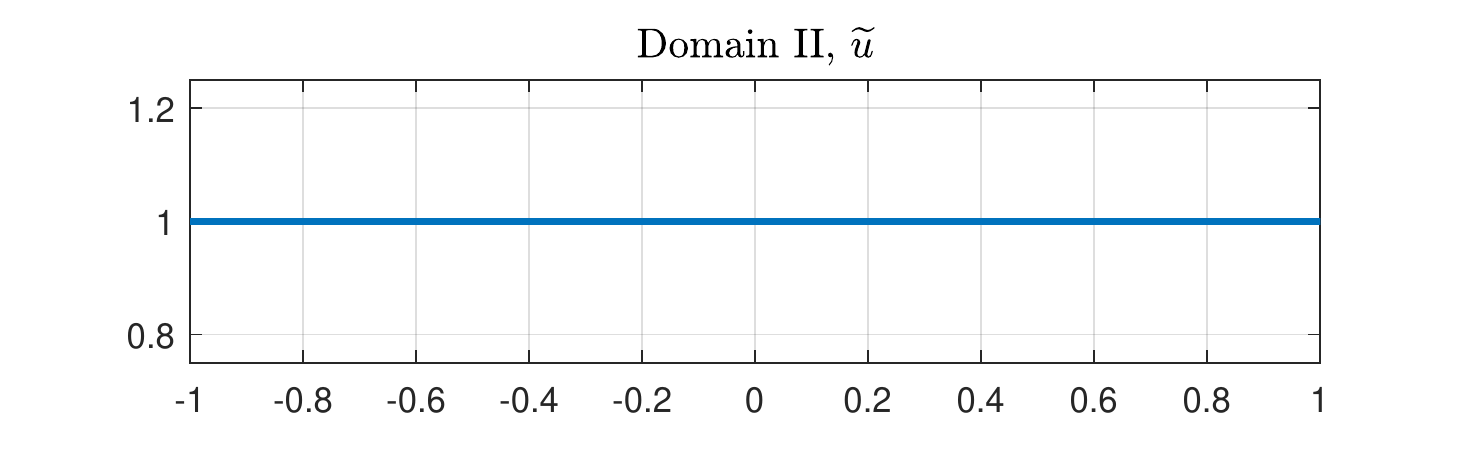}}
\hspace*{-0.25 cm}
\mbox{\includegraphics[width=0.525\textwidth]{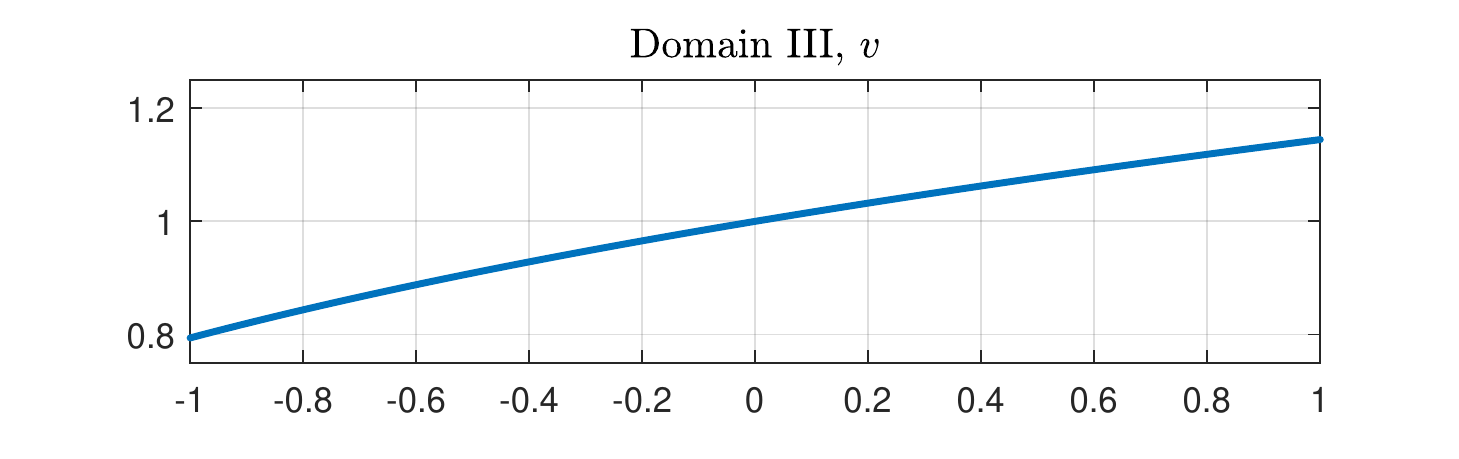}
\includegraphics[width=0.525\textwidth]{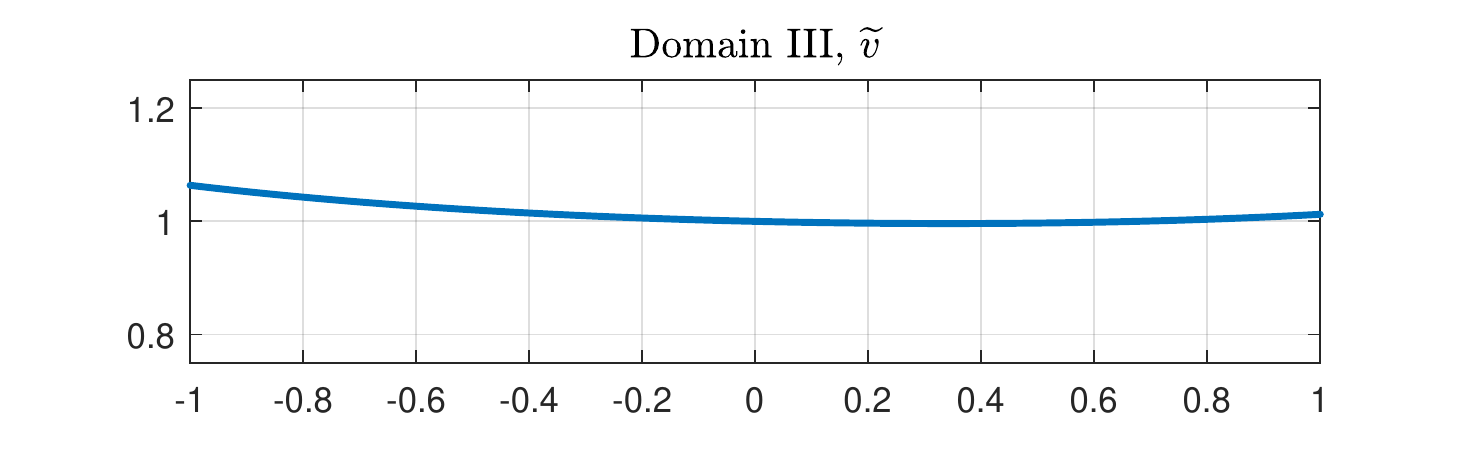}}
 \caption{The computed solutions to equations (\ref{hypergeom}), (\ref{hypergeomt1}), (\ref{II}), (\ref{infv2}) and (\ref{infv3}), all with the condition $y(0)=1$ and parameter values $a=-1/3$, $b = c = 1/2$,  the building blocks for the numerical construction of  the hypergeometric function on the real axis. Accuracy on the order of machine precision is achieved for each of these solutions with, respectively,  $n = 28$, $n = 28$, $n = 1$, $n = 30$ and $n = 27$ Chebychev coefficients. }
 \label{real_sols_fig}
\end{figure}

The solution to (\ref{infv2}) with $v(0)=1$ and $a=-1/3$, $b=c=1/2$ is 
$v=(1+s/2)^{1/3}$, see Figure~\ref{real_sols_fig}.  The 
solution $\widetilde{v}$ to (\ref{infv3}) with $\widetilde{v}(0)=1$, also plotted in Figure~\ref{real_sols_fig}, does not appear to have a 
simple closed form.  

\subsection{Matching at the domain boundaries}\label{matching}

In this section, we have so far shown (for the studied example) that we can compute solutions for each respective domain to essentially machine precision with about $n = 30$ Chebychev coefficients per domain. These solutions are analytic functions and are also the building blocks for  
the general solution to a Fuchsian equation, as per 
Frobenius (\ref{frobenius}). Note that in this 
approach infinity is just a normal point on the compactified real axis, thus 
large values of $|x|$ are not qualitatively different from points 
near $x=0$.

The construction of these analytic solutions allows us to continue the 
solution into domain I, which is just the 
hypergeometric function, to the whole real line, even to points where 
it is singular. This is done as follows: the general solution to the 
hypergeometric equation (\ref{hypergeom}) in domain II has the form
\begin{equation}
y_{II}(t)=\alpha u(t)+\beta t^{c-a-b} \widetilde{u}(t), \label{match1params}
\end{equation}
where $\alpha$ and $\beta$ are 
constants. These constants are determined by the condition that the 
hypergeometric function is differentiable at the boundary $x=1/2$ (which corresponds to $t = 1/2$ since $t = 1 - x$) 
between domains I and II:
\begin{equation}
y(1/2) = y_{II}(1/2), \qquad y'(1/2) = y'_{II}(1/2)\frac{\mathrm{d}t}{\mathrm{d}x} = - y'_{II}(1/2).  \label{match1}
\end{equation} 
The derivative of the numerical solutions at the endpoints of the interval can be computed using the formul\ae\ $T'_j(1) = j^2$ and $T'_j(-1) = (-1)^{j+1}j^2$. Alternatively, Chebfun can be used, which implements the recurrence relations in~\cite{MHCheb} for computing the derivative of a truncated Chebychev expansion. In the example studied, 
we find as expected $\alpha=0$ and $\beta=1$ up to a numerical error of $10^{-16}$.

\begin{remark}\label{remrank}
    If 
    $c-a-b\in \mathbb{Z}$ which is excluded by (\ref{generic}), it is 
    possible that the second solution is not linearly independent. This 
    would lead to non-unique values of $\alpha$ and $\beta$ in  (\ref{match1}). 
\end{remark}

\begin{figure}[htb!]
\hspace*{-0.35 cm}
\mbox{\includegraphics[width=0.525\textwidth]{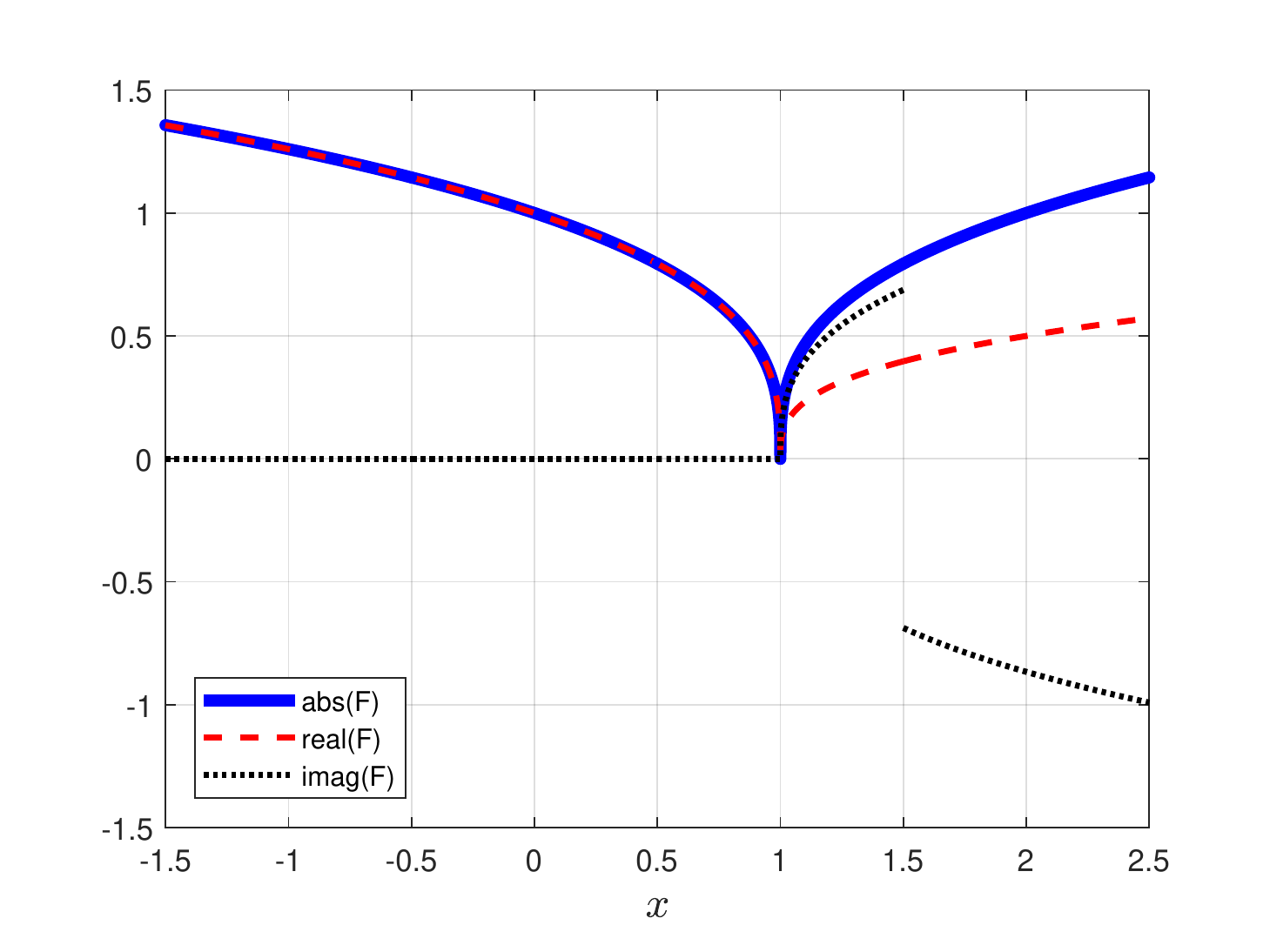}
\includegraphics[width=0.525\textwidth]{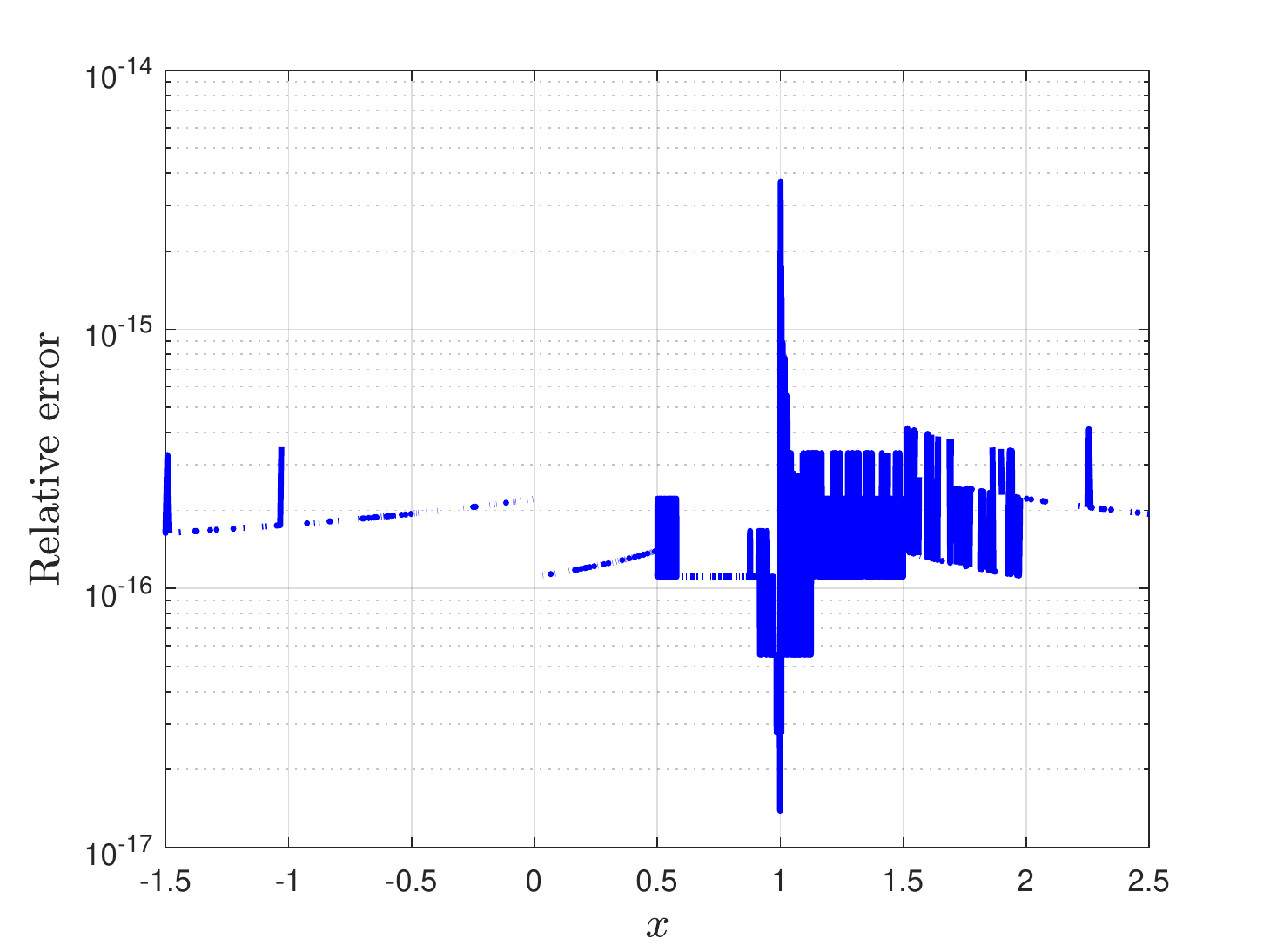}}
\hspace*{-0.35 cm}
\mbox{\includegraphics[width=0.525\textwidth]{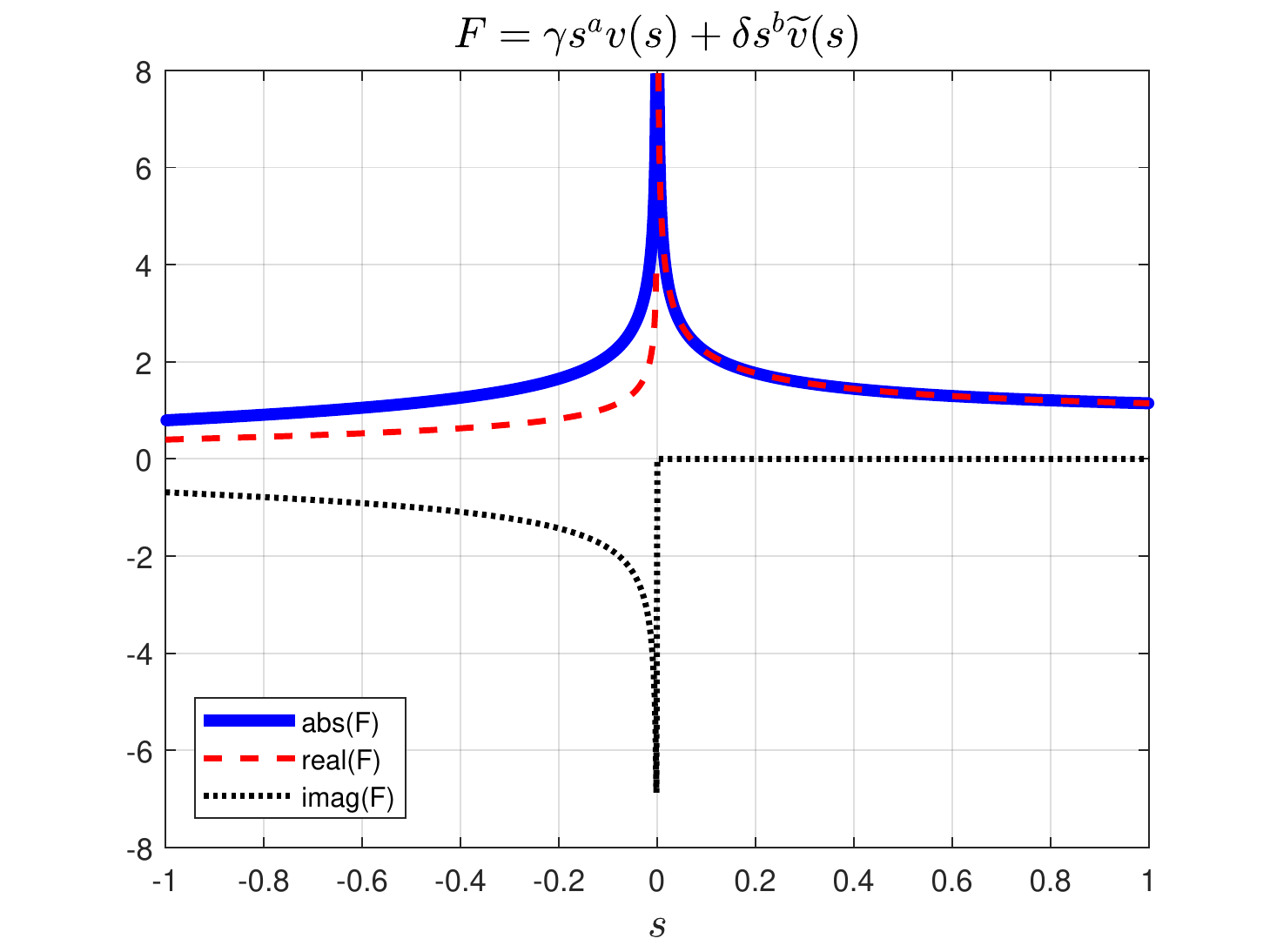}
\includegraphics[width=0.525\textwidth]{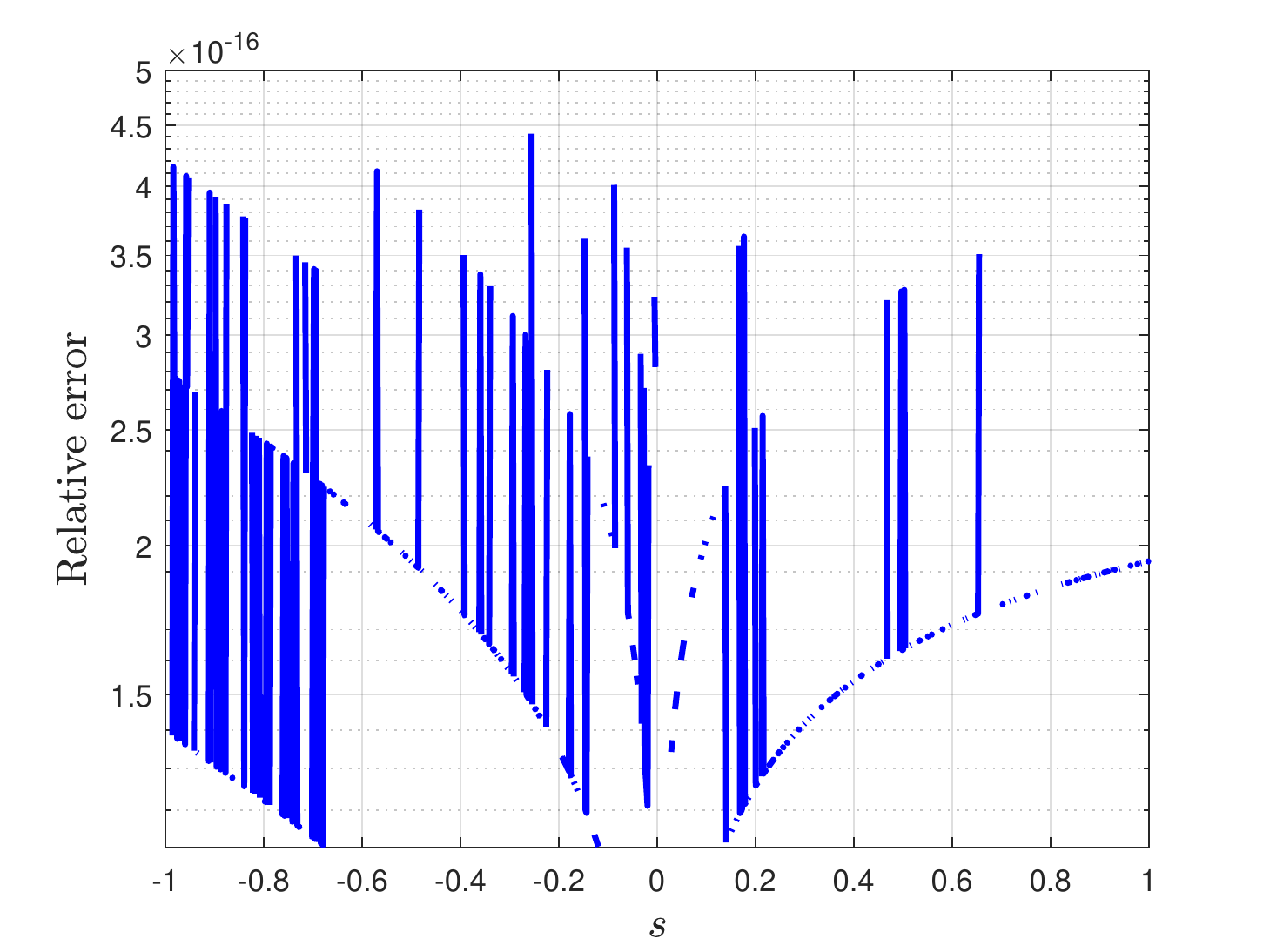}
}
 \caption{The hypergeometric function $F(-1/3, 1/2, 1/2, x)$  in the $x$ and $s$ planes (column 1), numerically constructed from the solutions in Figure~\ref{real_sols_fig},  and the relative error (column 2).}
 \label{realinterpol}
\end{figure}

In the same way the hypergeometric function can be analytically 
continued along the negative real axis. The general solution in 
domain III can be written as 
\begin{equation}
y_{III}=\gamma s^{a}v(s)+\delta 
s^{b}\widetilde{v}(s), \label{match2params}
\end{equation}
 with $\gamma$ and $\delta$ constants. 
Again linear independence of these solutions is assured by condition 
(\ref{generic}). The matching conditions at 
$x=-1/2$ (which corresponds to $s = 1$ since $s = -1/(x-1/2)$) are
\begin{equation}
y(-1/2) = y_{III}(1), \qquad y'(-1/2) = y'_{III}(1)\left[\frac{\mathrm{d}s}{\mathrm{d}x}\right]_{x=-1/2} =  y'_{III}(1).  \label{match2}
\end{equation} 
 For the 
studied example we find as expected $\gamma=1$ and $\delta=0$ with an 
accuracy better than $10^{-16}$.
Note that the hypergeometric function is in this way analytically 
continued also to positive values of $x\geq 3/2$, but this does not imply 
that the function is continuous at $x=3/2$ as can be seen in 
Figure~\ref{realinterpol}. The reason is the 
appearance of roots in the solutions, see Remark~\ref{rem1}, which 
leads to different branches of the hypergeometric function (Matlab chooses different  branches of the functions $t^{c-a-b} = t^{1/3}$ and $s^{a} = s^{-1/3}$ in domains II and III, respectively, causing the discontinuity in the imaginary part of the solution in Figure~\ref{realinterpol}). 


It can be seen in the top-right frame that full precision is attained on the interval except in the vicinity of $x=1$ since the singularity causes function evaluation to be ill-conditioned in this neighbourhood.  The second row of Figure~\ref{realinterpol} illustrates the computed hypergeometric function and the error in the $s$-plane (recall that $s \in [-1, 0]$ is mapped to $x \in [3/2, +\infty)$ and $s \in (0, 1]$ is mapped to $x \in (-\infty, -1/2])$). We have thus computed the hypergeometric function for the test example on the whole compactified real line to essentially machine precision. To recapitulate, this required the solution of five almost banded linear systems of the form shown in Figure~\ref{realI}, followed by the imposition of continuous differentiability at the domain boundaries as in (\ref{match1}) and (\ref{match2}).   

\section{Numerical construction of the hypergeometric function in the 
whole complex plane}\label{sec3}
\begin{figure}[htb!]
   \includegraphics[width=0.75\textwidth]{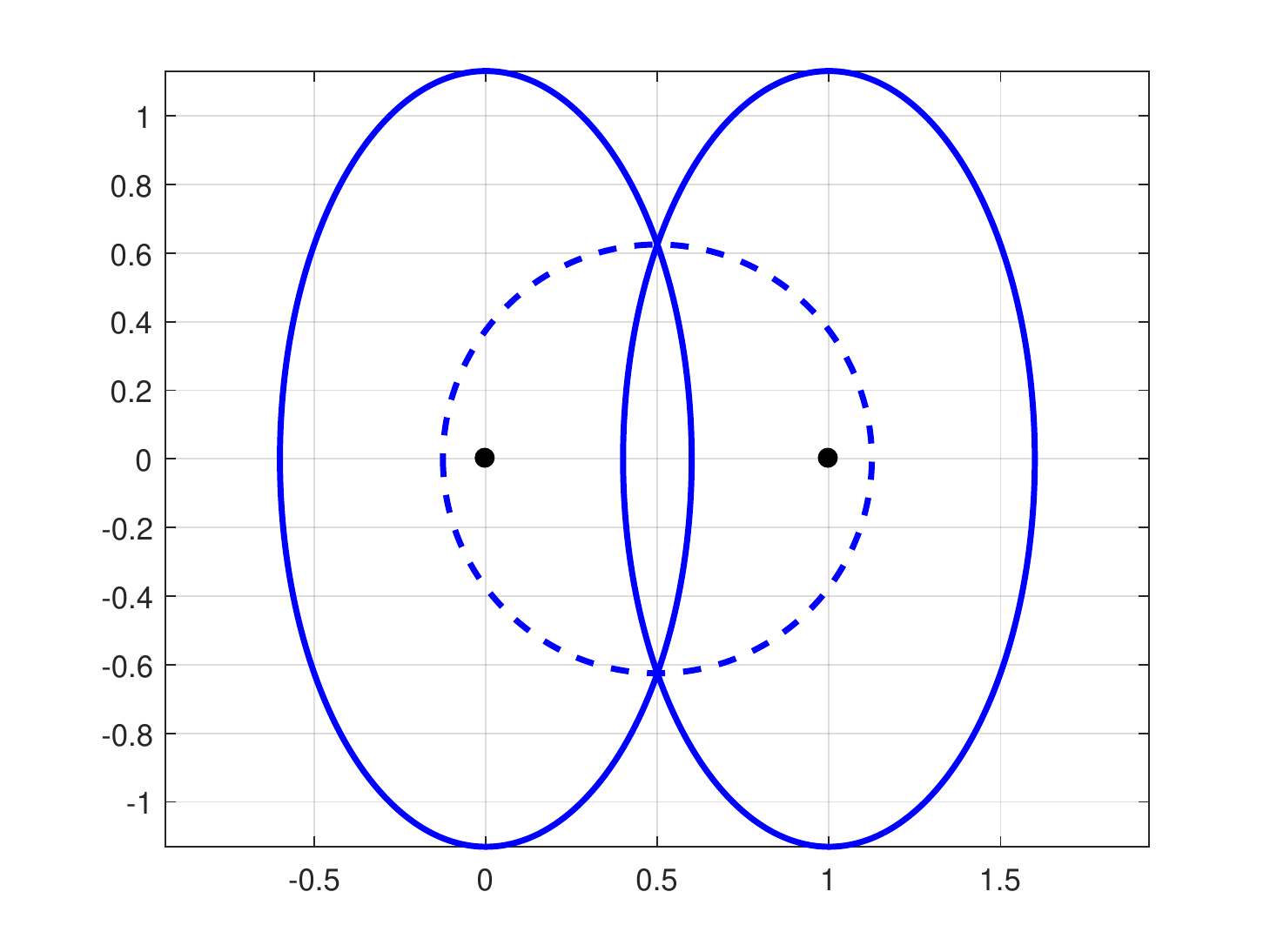}
 \caption{Domains for the computation of the hypergeometric function: 
 domain~I is the interior of the ellipse centered at \mbox{$z=0$}, 
 domain~II is the interior of the ellipse centered at \mbox{$z=1$}, and 
 domain~III is the exterior of the dashed circle centered at \mbox{$z=1/2$}.}
 \label{compldomains}
\end{figure}

In this section, the approach of the previous section is extended to 
the whole complex plane. Again three domains are introduced each of 
which exactly contains one of the three singular points 0, 1, and 
infinity, and which cover now the whole complex $z$ plane. To keep the 
number of domains to three and in order to have simply connected 
boundaries, we choose ellipses as shown in Figure~\ref{compldomains}:\\
- domain I: interior of the ellipse given by $\Re z =A\cos \phi$, 
$\Im z=B\sin\phi$, $\phi\in[-\pi,\pi]$.\\
- domain II: interior of the ellipse given by $\Re z=1+A\cos \phi$, 
$\Im z=B\sin\phi$, $\phi\in[-\pi,\pi]$.\\
- domain III: exterior of the circle $\Re z=1/2+R\cos \phi$, 
$\Im z=R\sin\phi$, $\phi\in[-\pi,\pi]$, where $R=B\sqrt{1-1/(4A^{2})}$.\\
Thus each of the ellipses is centered at one of the singular points. 
The goal is to stay away  from the other singular points 
since the equation to be solved is singular there which might lead to 
numerical problems if one gets too close. As discussed in 
\cite{RS2,RS3}, a distance of the order of $10^{-3}$ is not 
problematic with the used methods, but slightly larger distances can 
be handled with less resolution. We choose $A$ and $B$ such that the shortest distances between the boundaries of domains I, II and III and the singular points are equal. In the $s$-plane, the singular points $z = 0$ and $z = 1$ are mapped to, respectively, $s = 2$ and $s = -2$  and the domain boundary is a circle of radius $1/R$ centred at the origin, and thus we require, for a given $1/2 < A < 1$, that $B$ is chosen such that
\begin{equation}
 1 - A = 2 - \frac{1}{R}, \quad R = B\sqrt{1-1/(4A^{2})} \quad \Rightarrow \qquad B = \frac{1}{(A+1)\sqrt{1 - 1/(4A^{2})}}. \label{ellparams}
\end{equation}  
For example, in Figure~\ref{compldomains}, $A = 0.6$ (the parameter value we use throughout) and thus the shortest distance from any domain boundary to the nearest singularity is $0.4$.   This allows us to cover the whole complex  plane 
whilst staying clear of the singularities.  
There are parts of the complex plane covered by more than one domain, 
the important point is, however, that the whole plane is covered. 

The solution in each of the ellipses is then constructed in 3 
steps:
\begin{itemize}
    \item[i)]   The code for the real axis described in the previous section is run 
on larger domains than needed for a computation on the real axis only 
in order to obtain the boundary values of the five considered 
forms of (\ref{hypergeom}) at the intersections of the ellipses with the real 
axis. 

    \item[ii)]  On the ellipse, the equivalent forms of 
    (\ref{hypergeom}) of the previous section 
are solved in the considered domain with boundary values given on the 
real axis, again with the US spectral method.

    \item[iii)]  The obtained solutions on the ellipses 
serve as boundary values for the solution to the Laplace equation in 
the interior of the respective domains. In this way, the solutions on 
the real axis are analytically continued to the complex domains. As described below, the Laplace equation is solved by representing the solution in the Chebychev--Fourier basis, which reduces the problem to a coupled (on an ellipse) or uncoupled (on a disk) system of second-order boundary value problems (BVPs) which we again solve with the US method.

\end{itemize}
In the last step   the matching described in subsection 
\ref{matching} provides the hypergeometric function on the whole 
Riemann sphere built from the 
constructed holomorphic function in the three domains.
 As 
detailed below, this can be achieved as before with 
spectral accuracy as will be again discussed for the example 
$F(-1/3,1/2,1/2,z)$.

\subsection{Domain I}
In domain I, the task is to give the solution to equation 
(\ref{hypergeom}) with $y(0)=1$. In step i) the 
solution is first constructed on the interval $x\in[-A,A]$ which 
yields $F(a,b,c,A)$ with the US method detailed in the previous section. We find that $n = 40$ is sufficient to compute the solution to machine precision. 

In step ii), the ODE (\ref{hypergeom}) with $x$ replaced by $z$ is solved on the ellipse 
\begin{equation}
    z(\phi) = \frac{A+B}{2}\exp(i\phi) + \frac{A-B}{2}\exp(-i\phi),\quad 
    \phi\in[-\pi,\pi]
    \label{zel}
\end{equation}
as an ODE in $\phi$,
\begin{equation}
    z_{\phi}z(1-z)y_{\phi\phi}+\left[z^2(1-z)+(c-(1+a+b)z)z_{\phi}^2\right]y_{\phi}-abz_{\phi}^{3}y=0
    \label{cIel},
\end{equation}
where an index $\phi$ denotes the derivative with respect to $\phi$, 
and where $z$ is given by (\ref{zel}). We seek the solution to this ODE 
with the boundary conditions 
$y(\phi=-\pi)=y(\phi=\pi)=F(a,b,c,-A)$. Since the solution is periodic in $\phi$, it is natural to apply Fourier methods to solve (\ref{cIel}). The Fourier spectral method, which we briefly present, is entirely analogous to the US spectral method---indeed it served as the inspiration for the US spectral method---but simpler since it doesn't require a change of basis. Note that all the variable coefficients in (\ref{cIel}) are band-limited functions of the form $\sum_{j = -m}^{m}a_j e^{i j \phi}$ with $m = 3$. As in the US method, we require multiplication operators to represent the differential equation in coefficient space.  Hence,  suppose $a(\phi) = \sum_{j = -m}^{m}a_j e^{i  j \phi}$ and $y(\phi) = \sum_{j = -\infty}^{\infty}y_j e^{i j \phi}$, then
\begin{equation*}
a(\phi) y(\phi) = \sum_{j = -m}^{m}a_j e^{i j \phi } \sum_{j = -\infty}^{\infty}y_j e^{i j \phi} = \sum_{j = -\infty}^{\infty}c_j e^{i j \phi},
\end{equation*}
where
\begin{equation*}
c_j = \sum_{k = j-m}^{j+m} a_{j-k}y_j, \qquad j \in \mathbb{Z},
\end{equation*}
or
\begin{equation*}
\bm{c} = \mathcal{T}[a(\phi)]\bm{y}:=
\left(
{\renewcommand{\arraystretch}{1.0}
\begin{array}{c c c c c c c}
\ddots &      & \ddots &        &        &        &        \\
       &  a_m & \cdots & a_{-m} &        &        &        \\
       &      &   a_m  & \cdots & a_{-m} &        &        \\     
       &      &        & a_m    & \cdots & a_{-m} &         \\
       &      &        &        & \ddots &        & \ddots
\end{array}
}
  \right)
  \left(
{\renewcommand{\arraystretch}{1.0}
\begin{array}{c }
\vdots \\
y_{-1} \\
y_0    \\
y_1    \\
\vdots
\end{array}
}
  \right).
\end{equation*}
In the Fourier basis the differential operators are diagonal:
\begin{equation*}
\mathcal{D}_1 = 
i \left(
{\renewcommand{\arraystretch}{1.0}
\begin{array}{c c c c c}
 \ddots &    &     &    &         \\
        & -1 &     &    &          \\
        &    &  0  &    &         \\
        &    &     & 1  &          \\
        &    &     &    & \ddots
\end{array}
}
  \right), \qquad \mathcal{D}_2 =  \mathcal{D}_1^2,
\end{equation*}
and thus in coefficient space equation (\ref{cIel}) without the boundary conditions becomes $\mathcal{L}\bm{y} = \bm{0}$, where
\begin{equation*}
\mathcal{L}:= \mathcal{T}[a_2(\phi)]\mathcal{D}_2 + \mathcal{T}[a_1(\phi)]\mathcal{D}_1 + \mathcal{T}[a_0(\phi)],
\end{equation*}
and $a_2(\phi)$, $a_1(\phi)$ and $a_0(\phi)$ denote the variable coefficients in (\ref{cIel}). To find the $2n+1$ coefficients of the solution, $y_j$, $j = -n, \ldots, n$,  we need to truncate $\mathcal{L}$, for which we define the $(2n+1) \times \infty$ operator
\begin{equation*}
\mathcal{P}_{-n,n} = \left( \bm{0}, I_{-n,n}, \bm{0} \right).
\end{equation*}
The subscripts of the $(2n+1) \times (2n+1)$ identity matrix $I_{-n,n}$ indicate the indices of the vector on which it operates, e.g., $\mathcal{P}_{-n,n}\bm{y} = \left[y_{-n}, \ldots, y_n \right]^{\top}$. Then the system to be solved to approximate the solution of (\ref{cIel}) is
\begin{equation}
\left(
{\renewcommand{\arraystretch}{2.0}
\begin{array}{c c c c}
(-1)^{-n} & (-1)^{1-n} & \cdots & (-1)^{n} \\
\multicolumn{4}{c}{\mathcal{P}_{-n,n-1} \mathcal{L} \mathcal{P}_{-n,n}^{\top}}
\end{array}}
\right)
\left(
\begin{array}{c}
y_{-n} \\
\vdots \\
y_n
\end{array}
\right)
=
\left(
\begin{array}{c}
F(a,b,c,-A)  \\
0            \\
\vdots       \\
0
\end{array}
\right)
.  \label{Fsyst}
\end{equation}

\begin{figure}[htb!]
\hspace*{-0 cm}
\mbox{
\includegraphics[width=0.5\textwidth]{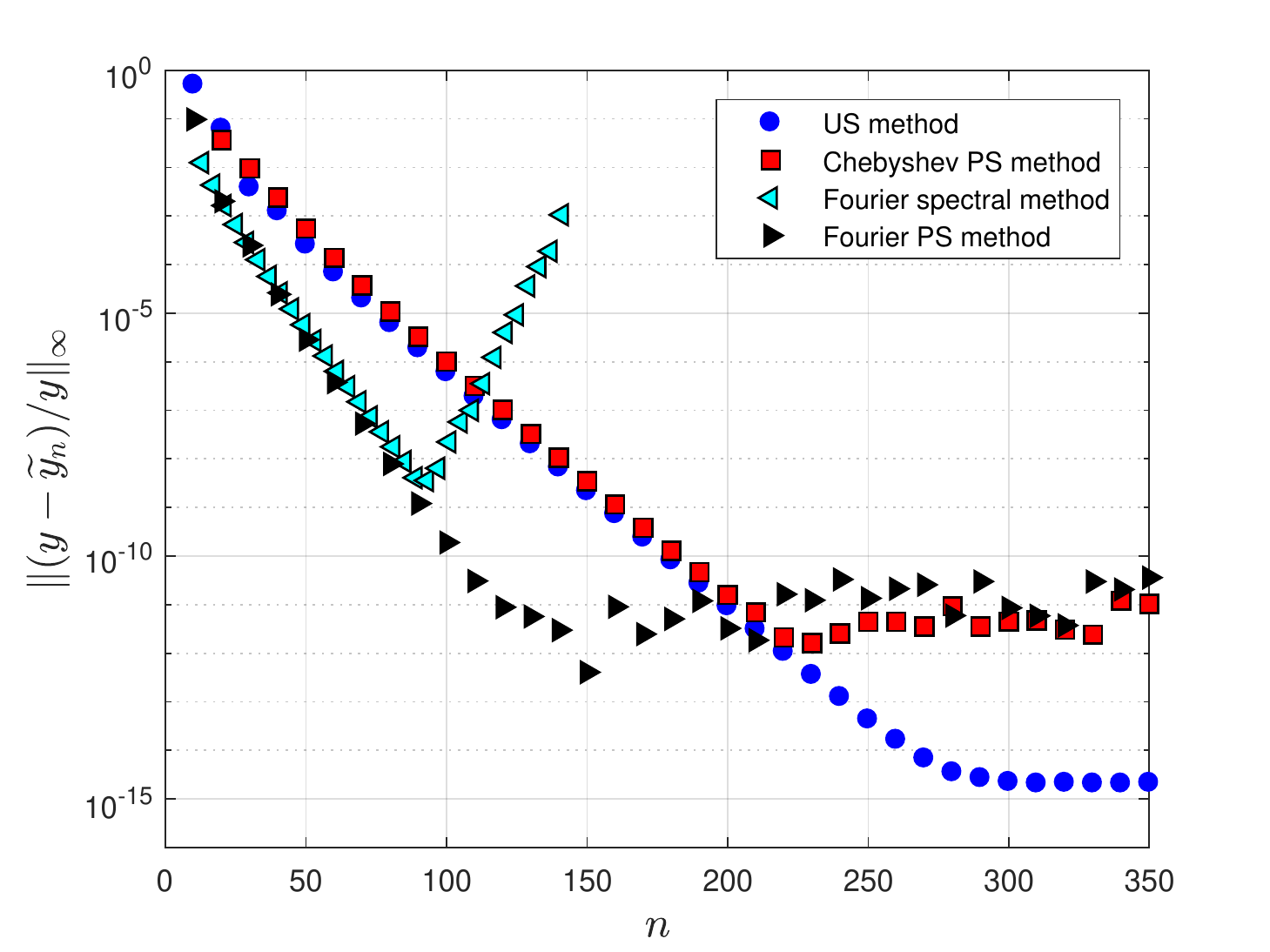}\includegraphics[width=0.5\textwidth]{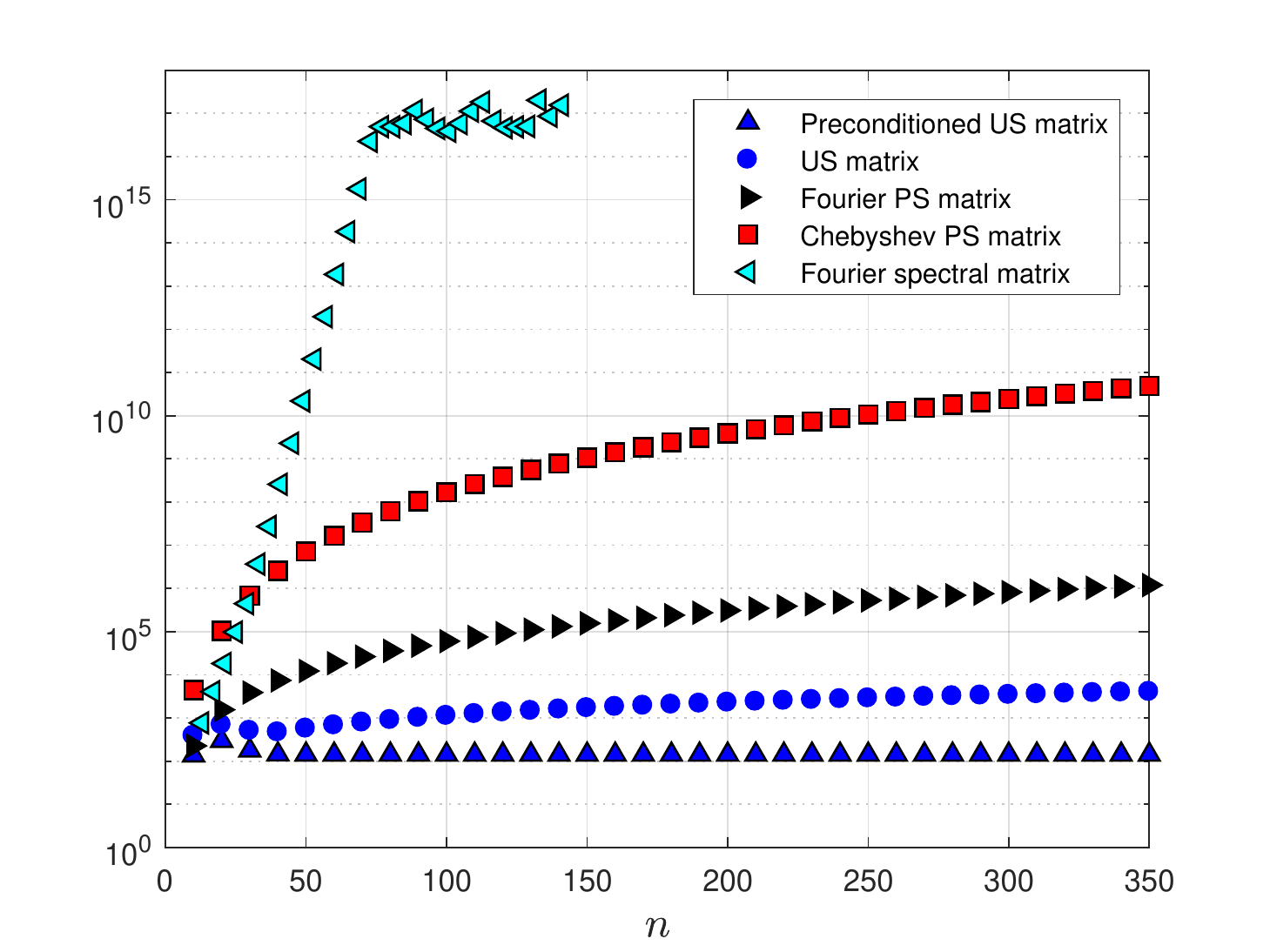}
}

 \caption{Numerical solution of the hypergeometric function 
 $F(-1/3,1/2,1/2,x)$ on the ellipse (\ref{zel}); on the left the 
 maximum relative error on the ellipse decreases exponentially with $n$, and on the right the growth of condition numbers of the matrices of the preconditioned US method, the US method, the Fourier PS method, the Chebychev PS method and the Fourier spectral method are, respectively, bounded for all $n$, linear, $\mathcal{O}\left(n^{2.4} \right)$, $\mathcal{O}\left(n^{4.5}\right)$  and $\mathcal{O}\left(e^{0.56 N} \right)$, where $N = 2n+1$.}
 \label{elI}
\end{figure}

In Figure~\ref{elI}, (\ref{cIel}) is also solved with a Fourier pseudospectral (PS) method~\cite{WR}, which operates on solution values at equally spaced points on $\phi \in [-\pi, \pi)$. The Fourier and Chebyshev PS methods used in Figure~\ref{elI} lead to dense matrices whereas the Fourier and Chebyshev spectral methods give rise to almost banded linear systems with bandwidths 3 and 35, respectively.  The almost banded Fourier spectral matrices (\ref{Fsyst}) have a single dense top row whereas the US matrices have two dense top rows (one row for each of the conditions $y(\phi=-\pi)=y(\phi=\pi)=F(a,b,c,A)$). 

Note that the Fourier methods converge at a faster rate than the Chebychev methods. This is to be expected since generally for periodic functions Fourier series converge faster than Chebychev series by a factor of $\pi/2$~\cite{trefweid}. Again the US method achieves the best accuracy and, as before, this is due to the difference in the conditioning of the methods, as shown in the right frame of Figure~\ref{elI}.

Unlike the equations solved in section~\ref{sec2},  (\ref{cIel}) has no singular points on its domain, and thus the condition numbers of the Chebyshev PS and US matrices grow at different rates than in Figure~\ref{realI}.  We find that, as shown in~\cite{US}, the condition numbers of the US matrices grow linearly with $n$ and the preconditioned US matrices have condition numbers that are bounded for all $n$. The condition numbers of the Fourier PS, Chebychev PS and Fourier spectral matrices grow as, respectively, $\mathcal{O}\left(n^{2.4} \right)$, $\mathcal{O}\left(n^{4.5}\right)$  and $\mathcal{O}\left(e^{0.56 N} \right)$ (where $N = 2n+1$, the number of Fourier coefficients of the solution in (\ref{Fsyst})\footnote{In Figure~\ref{elI}, the results for the Fourier spectral method are plotted against $N$ and not $n$, as the axis label indicates.
}),  according to a least squares fit of the data.  Figure~\ref{elI} shows that the exponential ill-conditioning of the Fourier spectral matrices results in a rapid loss of accuracy for large enough $N$. 

In step iii), in order to analytically continue the hypergeometric function to the interior 
of the ellipse, we use the fact that the function is holomorphic 
there and thus harmonic. Therefore we can simply solve the Laplace 
equation in elliptic coordinates
\begin{equation}
        z(r,\phi) = \frac{A+B}{2}r\exp(i\phi) + \frac{A-B}{2}r\exp(-i\phi),\quad 
   r\in[0,1],\quad \phi\in[0,2\pi]
    \label{zrphi}.
\end{equation}
In these coordinates, the Laplace operator reads
\begin{equation}
    \begin{split}
   r^2 \Delta&=\left(\frac{1}{A^{2}}-\frac{1}{B^{2}}\right)\left[\frac{1}{2}\cos 2\phi\left(r^2\partial_{rr}-r \partial_{r}
    -\partial_{\phi\phi}\right)+\sin2\phi\left(\partial_{\phi}-r\partial_{r\phi}\right)\right]\\
    &+\frac{1}{2}\left(\frac{1}{A^{2}}+\frac{1}{B^{2}}\right)\left(r^2\partial_{rr}+r\partial_{r}
    + \partial_{\phi\phi}\right).
    \end{split}
    \label{Lapel}
\end{equation}

Notice that (\ref{Lapel}) simplifies considerably on the disk (if $A = B$), which results in a more efficient numerical method. However, using ellipses allows us to increase the distance between the domain boundaries and the nearest singularities and we have found that, if the closest singularity is sufficiently strong, this yields more accurate solutions compared to using disks. For the test problem (\ref{example}), where the exponent of the singularity at $z = 1$ is $c - a - b = 1/3$, we have found that using an ellipse as opposed to a disk improves the accuracy only by a factor slightly more than two. However, for an example to be considered in section~\ref{examples_sect}  (the first three rows of Table~\ref{table2}) where the exponent at $z = 1$ is $c - a -b = -0.6$, using an ellipse (with $A = 0.6$ in (\ref{ellparams})) yields a solution that is more accurate than the solution obtained on a disk (with parameters $A = B = 0.7574\ldots$, obtained by solving (\ref{ellparams}) with $A = B$) by two orders of magnitude. 

Another possibility, which combines the advantages of ellipses (better accuracy) and disks (more efficient solution of the Laplace equation), is to conformally map disks to ellipses as in~\cite{atkinsonellipse}. However, we found that computing this map (which involves elliptic integrals) to machine precision for $A = 0.6$ in (\ref{ellparams}) requires more than $1200$ Chebyshev coefficients. This is about four times the number of Chebyshev coefficients required to resolve the solution in Figure~\ref{elI}. In addition, the first and second derivatives of the conformal map, which are needed to solve the hypergeometric equation (\ref{hypergeom}) and also (\ref{hypergeomt1})--(\ref{II}) on ellipses, involve square roots and this requires that the right branches be chosen. Hence, due to the expense and complication of this approach we did not pursue it further.
%

Yet another alternative is to use rectangular domains, where the boundary data have to be generated by 
solving ODEs on the 4 sides of each rectangle. Then the solution can be expressed as a bivariate Chebychev expansion. A disadvantage of this approach, noted in~\cite{boydcomp}, is that the grid clusters at the four corners of the domain which decreases the efficiency of the method. 

To obtain the numerical solution of the Laplace equation on the ellipses, we use the ideas behind the optimal complexity Fourier--Ultraspherical spectral method in~\cite{Weather} for the disk. Since the solution is periodic in the angular variable $\phi$, it is approximated by a radially dependent truncated Fourier expansion:
\begin{equation}
y \approx \widehat{y}_m(r,\phi) := \sum_{k = -m/2}^{m/2 -1} u_k(r) e^{ik\phi}, \qquad r \in  [-1, 1], \qquad \phi \in [-\pi, \pi).  \label{laplacesol}
\end{equation} 
As suggested in~\cite{Weather,trefethen}, we let $r \in [-1, 1]$ instead of $r \in [0, 1]$ to avoid excessive clustering of points on the Chebychev--Fourier grid near $r = 0$. With this approach the origin $r = 0$ is not treated as a boundary. Since $(r, \phi)$ and $(-r, \phi + \pi)$ are mapped to the same points on the ellipse, we require that 
\begin{equation}
\widehat{y}_m(-r,\phi + \pi) = \widehat{y}_m(r,\phi). \label{bmcII}
\end{equation}
On the boundary of the ellipse we specify $\widehat{y}_m(1,\phi) = \widetilde{y}_n(\phi)$, where $\widetilde{y}_n(\phi)$ is the approximate solution of (\ref{cIel}) obtained with the US method. Suppose that $\widetilde{y}_n(\phi)$ has the Fourier expansion $\widetilde{y}_n(\phi) = \sum_{k = -\infty}^{\infty} \gamma_k e^{ik\phi}$. Using the property (\ref{bmcII}), the boundary condition $\widehat{y}_m(1,\phi) = \widetilde{y}_n(\phi)$ becomes
\begin{equation}
u_k(1) = \gamma_k, \qquad u_k(-1) = (-1)^k\gamma_k, \qquad k = -m/2, \ldots, m/2 - 1.  \label{LaplBCs}
\end{equation}
Substituting (\ref{laplacesol}) into (\ref{Lapel}), we find that the Laplace equation $r^2 \Delta \widehat{y}_m = 0$ reduces to the following coupled system of BVPs, with boundary conditions given by  (\ref{LaplBCs}):
\begin{equation}  \label{BVPsyst}
\begin{split}
&\frac{1}{4}\left( \frac{1}{A^2} - \frac{1}{B^2} \right)\left\lbrace r^2 u''_{k-2} - \left[ 1 + 2(k-2) \right] r u'_{k-2} + (k-2)\left[(k-2) + 2\right] u_{k-2} \right\rbrace  \\
& + \frac{1}{2}\left( \frac{1}{A^2} + \frac{1}{B^2} \right)\left\lbrace r^2 u''_{k} + r u'_{k} - k^2 u_{k} \right\rbrace + \\
&  \frac{1}{4}\left( \frac{1}{A^2} - \frac{1}{B^2} \right)\left\lbrace r^2 u''_{k+2} - \left[ 1 + 2(k+2) \right] r u'_{k+2} + (k+2)\left[(k+2) - 2\right] u_{k+2} \right\rbrace = 0,
\end{split}
\end{equation}
for $k = -m/2, \ldots, m/2 -1$, where  $u_{k} = 0$ if $k < -m/2$ or $k > m/2 - 1$. Note that on a  disk ($A = B$) the system (\ref{BVPsyst}) reduces to a decoupled system of $m$ BVPs. The BVPs are solved using the US method, as in section~\ref{sec2}. Let 
\begin{equation}
\mathcal{T}_0 = \mathcal{S}_1\mathcal{S}_0, \qquad 
\mathcal{T}_1 = \mathcal{S}_1\mathcal{M}_1[r]\mathcal{D}_1, \qquad
\mathcal{T}_2 = \mathcal{M}_2[r^2]\mathcal{D}_2,  \label{termsdef}
\end{equation}
where the operators  in (\ref{termsdef}) are defined in section~\ref{sec2}. Let $\bm{u}^{(k)}$ denote the infinite vector of Chebychev coefficients of $u_k$,  then in coefficient space (\ref{BVPsyst}) becomes 
\begin{equation}  \label{BVPsystop}
\begin{split}
&\underbrace{\frac{1}{4}\left( \frac{1}{A^2} - \frac{1}{B^2} \right)\left\lbrace \mathcal{T}_2 - \left[ 1 + 2(k-2) \right] \mathcal{T}_1 + (k-2)\left[(k-2) + 2\right] \mathcal{T}_0  \right\rbrace}_{=\mathcal{L}^{(k-2)}} \bm{u}^{(k-2)} \\
& + \underbrace{\frac{1}{2}\left( \frac{1}{A^2} + \frac{1}{B^2} \right)\left\lbrace \mathcal{T}_2 + \mathcal{T}_1 - k^2\mathcal{T}_0  \right\rbrace }_{=\mathcal{M}^{(k)}}\bm{u}^{(k)}  + \\
&  \underbrace{\frac{1}{4}\left( \frac{1}{A^2} - \frac{1}{B^2} \right)\left\lbrace \mathcal{T}_2 - \left[ 1 + 2(k+2) \right] \mathcal{T}_1 + (k+2)\left[(k+2) - 2\right] \mathcal{T}_0  \right\rbrace}_{=\mathcal{R}^{(k+2)}} \bm{u}^{(k+2)}  = \bm{0}.
\end{split}
\end{equation}
The operators $\mathcal{L}^{(k-2)}$, $\mathcal{M}^{(k)}$ and $\mathcal{R}^{(k+2)}$ defined (\ref{BVPsystop}) are truncated and the boundary conditions (\ref{LaplBCs}) are imposed as follows to obtain a linear system for the first $n$ Chebychev coefficients of $u_k$, i.e., $\mathcal{P}_{n}\bm{u}^{(k)}$, for $k= -m/2, \ldots, m/2 - 1$:
\begin{equation}
\left(
{\renewcommand{\arraystretch}{4}
\begin{array}{c | c | c}
L^{(k-2)}_n & M^{(k)}_n & R^{(k+2)}_n 
\end{array}
}
\right)
\left(
{\renewcommand{\arraystretch}{2}
\begin{array}{c}
\mathcal{P}_{n} \bm{u}^{(k-2)} \\
\hline
\mathcal{P}_{n} \bm{u}^{(k)} \\
\hline
\mathcal{P}_{n} \bm{u}^{(k+2)}
\end{array}}
\right)
= \left(
\begin{array}{c}
\gamma_k \\
(-1)^k\gamma_k\\
0 \\
\vdots \\
0
\end{array}
\right),   \label{matrixeqs}
\end{equation} 
where
\begin{equation*}
L^{(k-2)}_n = \left(
\begin{array}{c c c c}
0 & \cdots & \cdots & 0 \\
0 & \cdots & \cdots & 0 \\
\multicolumn{4}{c}{\mathcal{P}_{n-2} \mathcal{L}^{(k-2)} \mathcal{P}_{n}^{\top}}
\end{array}
\right), \qquad
R^{(k+2)}_n = 
\left(
\begin{array}{c c c c}
0 & \cdots &  \cdots & 0 \\
0 & \cdots & \cdots & 0 \\
\multicolumn{4}{c}{\mathcal{P}_{n-2} \mathcal{R}^{(k+2)} \mathcal{P}_{n}^{\top}}
\end{array}
\right),
\end{equation*}
and
\begin{equation*}
M^{(k)}_n = \left(
\begin{array}{c c c c}
T_0(1) & T_1(1) & \cdots & T_{n-1}(1) \\
T_0(-1) & T_1(-1) & \cdots & T_{n-1}(-1)\\
\multicolumn{4}{c}{\mathcal{P}_{n-2} \mathcal{M}^{(k)} \mathcal{P}_{n}^{\top}}
\end{array}
\right).
\end{equation*}
The equations (\ref{matrixeqs}) can be assembled into two $nm/2 \times nm/2$ block tridiagonal linear systems: one for even $k$ and another for odd $k$ (recall that $\mathcal{P}_n \bm{u}^{(k)} = \bm{0}$ for $k < -m/2$ or $k > m/2-1$).  The systems can be further reduced by a factor of $2$ by using the fact that the function $u_k(r)$ has the same parity as $k$~\cite{Weather} because of the property (\ref{bmcII}). That is, if $k$ is even/odd, then $u_k$ is an even/odd function and hence only the even/odd-indexed Chebychev coefficients of $u_k$ are nonzero (and thus one of the two top rows imposing the boundary conditions in (\ref{matrixeqs}) may also be omitted) . Then the equations (\ref{matrixeqs}) are reduced to two $nm/4 \times nm/4$ block tridiagonal linear systems in which each off-diagonal block is tridiagonal and the diagonal block is almost banded with bandwidth one and a single dense top row.  On a disk, e.g., on domain~III, only the diagonal blocks of the system remain and the equations reduce to $m$ times $n/2 \times n/2$ tridiagonal plus rank one systems, which can be solved in $\mathcal{O}(n)$ operations with the Sherman-Morrison formula~\cite{Weather} resulting in a total computational complexity of $\mathcal{O}\left( m n \right)$. 

\begin{figure}[htb!]
\includegraphics[width=0.75\textwidth]{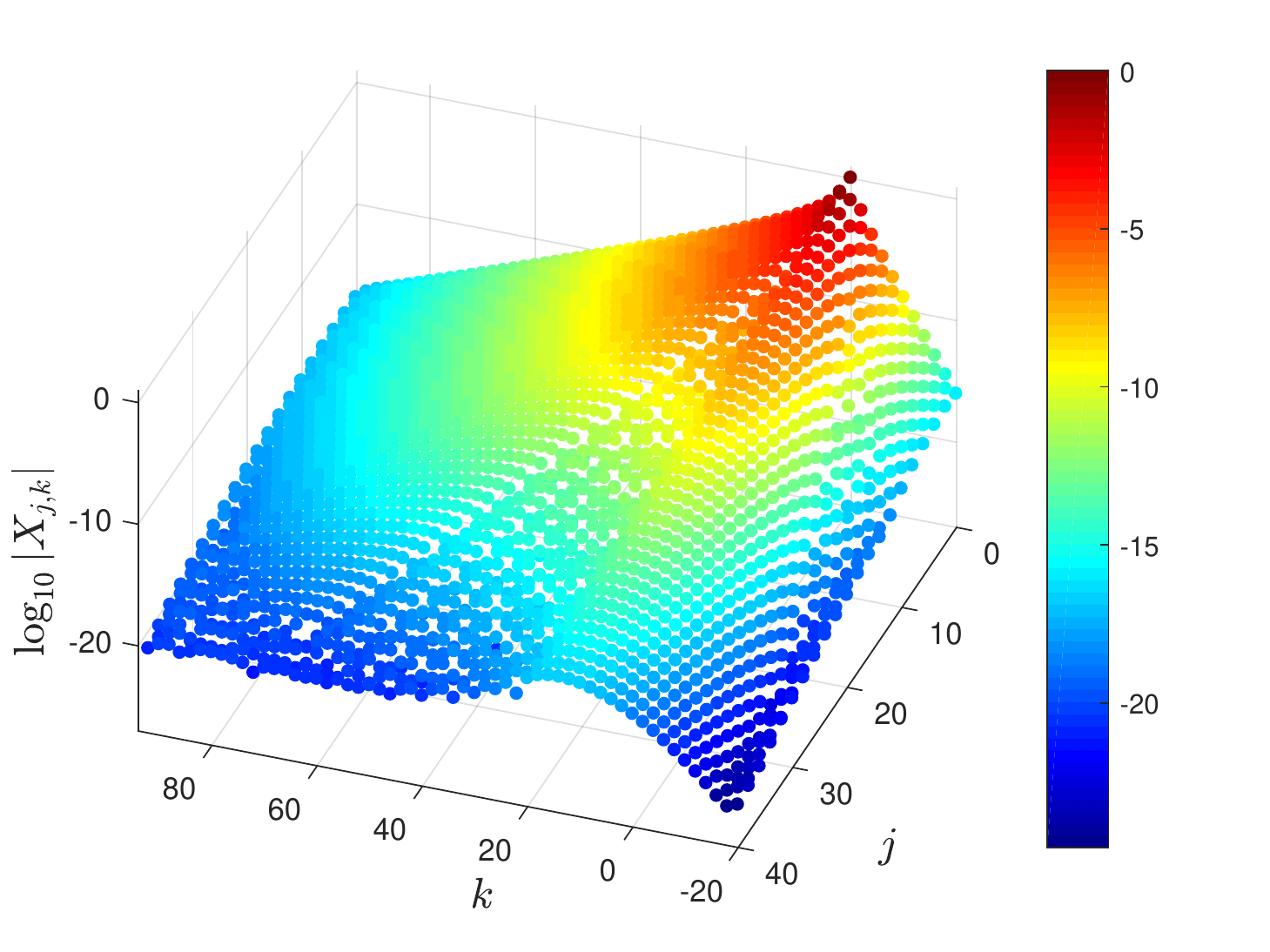}
 \caption{The magnitude of the Chebychev--Fourier coefficients of the hypergeometric function $F(-1/3, 1/2, 1/2, z)$, obtained by solving the Laplace equation on domain I.}
 \label{CFcoeffs}
\end{figure}

Solving the above system, the first $n$ Chebychev coefficients of $u_k$, where $k = -m/2, \ldots, m/2-1$, are obtained which are stored in column $k$ of an $n \times m$ matrix $X$ of Chebychev--Fourier coefficients. Then the solution expansion (\ref{laplacesol}) is approximated by 
 \begin{equation}
\widehat{y}_m(r,\phi) \approx  \sum_{k = -m/2}^{m/2 -1}\sum_{j = 0}^{n-1} X_{j,k} T_j(r) e^{ik\phi}, \qquad r \in  [-1, 1], \qquad \phi \in [-\pi, \pi). \label{CFexp} 
\end{equation} 
Figure~\ref{CFcoeffs} shows the exponential decrease in the magnitude of $X_{j,k}$ for the solution on domain I of the test problem (\ref{example}) with $m = 188$. Notice that $k$ ranges over only $k = -20, \ldots 93 = m/2-1$ instead of $k = -m/2, \ldots, m/2 - 1$ since we only set up the systems (\ref{matrixeqs}) for $k$ such that $\vert \gamma_k \vert $ is above machine precision. We use Chebfun (which uses the Fast Fourier Transform (FFT)) to compute the Fourier coefficients $\gamma_k$ of the function $\widetilde{y}_n$ on the domain boundary obtained in step (ii).  

To evaluate the Chebychev--Fourier expansion (\ref{CFexp}) at the set of $n_r n_{\phi}$ points $(r_i, \phi_j)$, $i = 1, \ldots, n_r$, $j = 1, \ldots, n_{\phi}$, where $0 \leq r_i \leq 1$, $ -\pi \leq \phi_j < \pi $, we form the $n_r\times 1$ and $n_{\phi} \times 1$ vectors $\bm{r}$ and $\bm{\phi}$ and compute the $n_r \times n_{\phi}$ matrix 
\begin{equation*}
\left(
{\renewcommand{\arraystretch}{4}
\begin{array}{c | c | c | c}
T_0(\bm{r}) & T_1(\bm{r}) & \cdots & T_{n-1}(\bm{r}) 
\end{array}
}
\right) X \left(
{\renewcommand{\arraystretch}{1.25}
\begin{array}{c}
\exp(-m/2\bm{\phi^{\top}}) \\
\hline
\exp((-m/2+1)\bm{\phi^{\top}}) \\
\hline
\vdots  \\
\hline
\exp((m/2-1)\bm{\phi^{\top}})
\end{array}}
\right).
\end{equation*}
The columns $T_j(\bm{r})$ are computed using the three term recurrence relation $T_{j+1} = 2rT_j - T_{j-1}$, with $T_0 =1 $ and $T_1(r) = r$. Alternatively, the expansion can be evaluated using barycentric interpolation~\cite{barycentric} in both the Chebychev and Fourier bases (which also requires the Discrete Cosine Transform and FFT to convert the Chebychev--Fourier coefficients to values on the Chebychev--Fourier grid) or by using Clenshaw's algorithm in the Chebychev basis~\cite{MHCheb} and Horner's method in the Fourier basis.

\begin{figure}[htb!]
\includegraphics[width=0.75\textwidth]{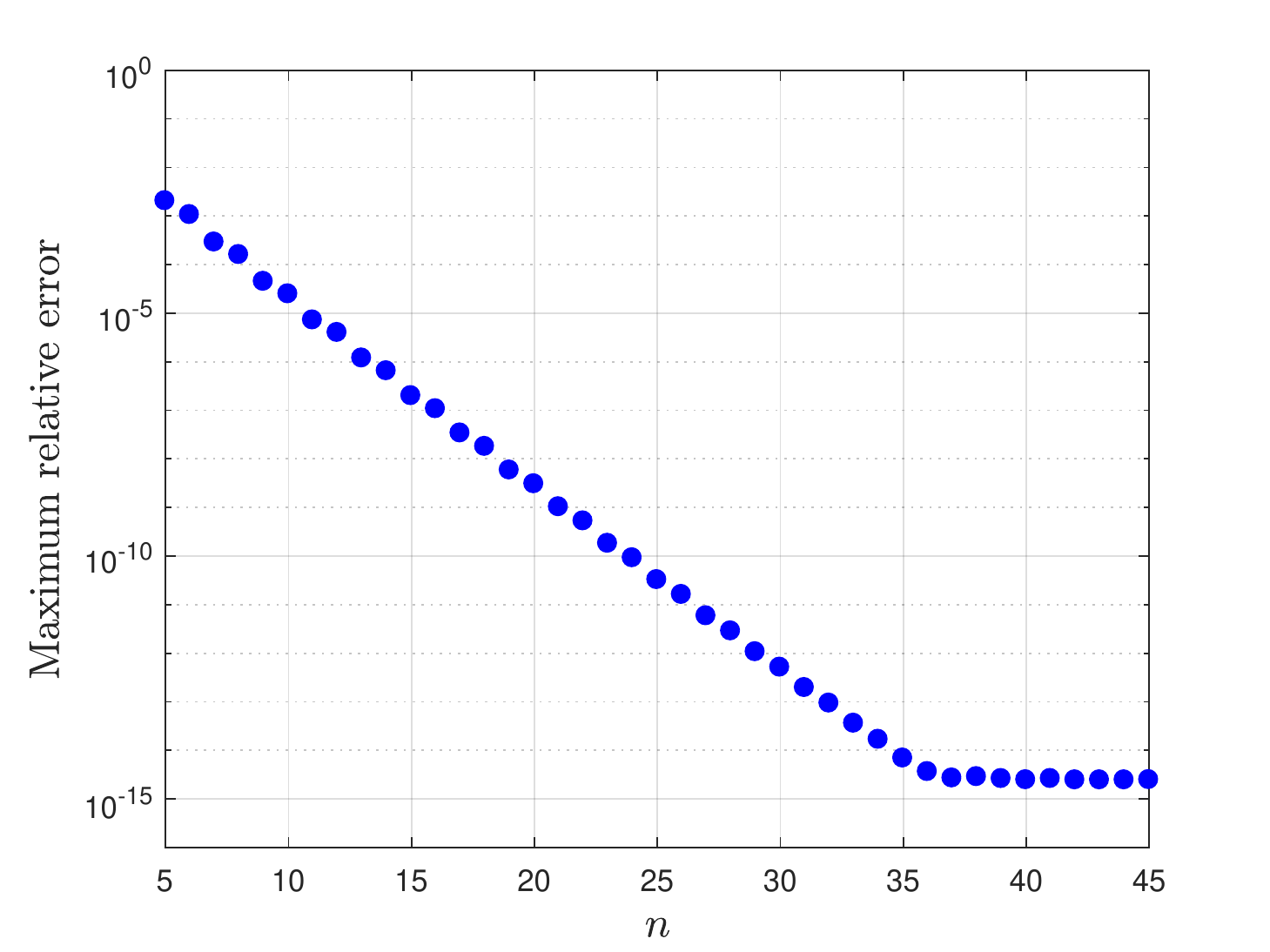}
 \caption{Spectral convergence of the Ultraspherical--Fourier spectral method to $F(-1/3, 1/2, 1/2, z)$ on domain I.}
 \label{elI2}
\end{figure}

Figure~\ref{elI2} shows the maximum relative error on domain I,  as measured on a $500 \times 500$ equispaced  grid on $(r, \phi) \in [0, 1]\times [-\pi, \pi)$,  as a function of $n$, the number of Chebychev coefficients of $u_k(r)$, $k = -m/2, \ldots, m/2-1$, for $m = 188$. 

\subsection{Domains II and III}

For the remaining domains and equations, the approach is the same: equations (\ref{hypergeomt1}) and (\ref{II}) are first solved on $[-A, A]$ (with $u(0) = 1 = \widetilde{u}(0)$), then on the ellipse centred at $z = 1$ shown in Figure~\ref{compldomains} and finally the Laplace equation is solved twice on the same ellipse but with different boundary data. Equations (\ref{infv2}) and (\ref{infv3}) are first solved on $[-1/R, 1/R]$ (with $v(0) = 1 = \widetilde{v}(0)$), then on the disk centred at $z=1/2$ shown in Figure~\ref{compldomains}, which is mapped to a disk of radius $1/R$ in the $s$-plane, and finally the Laplace equation is solved twice on a disk in the $s$-plane with different boundary data. The results are very similar to those obtained in Figures~\ref{elI} and~\ref{elI2}.

Since the solutions constructed in the present section are just 
analytic continuations of the ones on the real axis of the previous 
section, the hypergeometric function is built from it as in 
subsection \ref{matching}. Even the values of $\alpha$, $\beta$ in 
(\ref{match1params})--(\ref{match1}) and $\gamma$, $\delta$ in (\ref{match2params})--(\ref{match2}) are the same 
and can be 
taken from the computation on the real axis. Thus we have 
obtained the hypergeometric function in the three domains of 
Figure~\ref{compldomains} which cover the whole Riemann sphere. The 
computational cost in constructing it is essentially given by 
inverting five times the matrix approximating the Laplace operator 
(\ref{Lapel}), which can be performed in parallel (the one-dimensional computations are in comparison for 
free).

\begin{figure}[htb!]
     \hspace*{-0 cm}
 \mbox{   
 \includegraphics[width=0.5\textwidth]{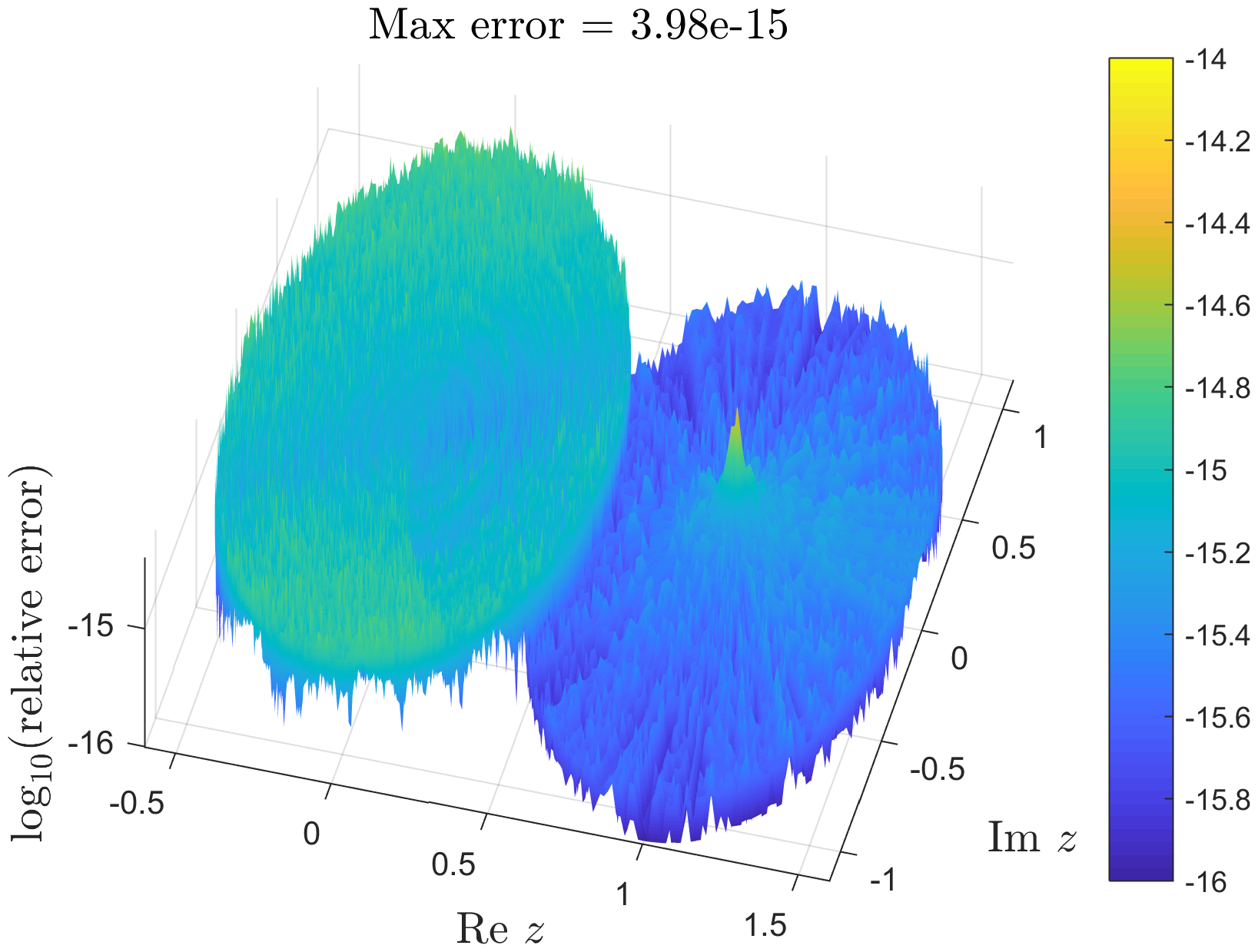}
 \includegraphics[width=0.5\textwidth]{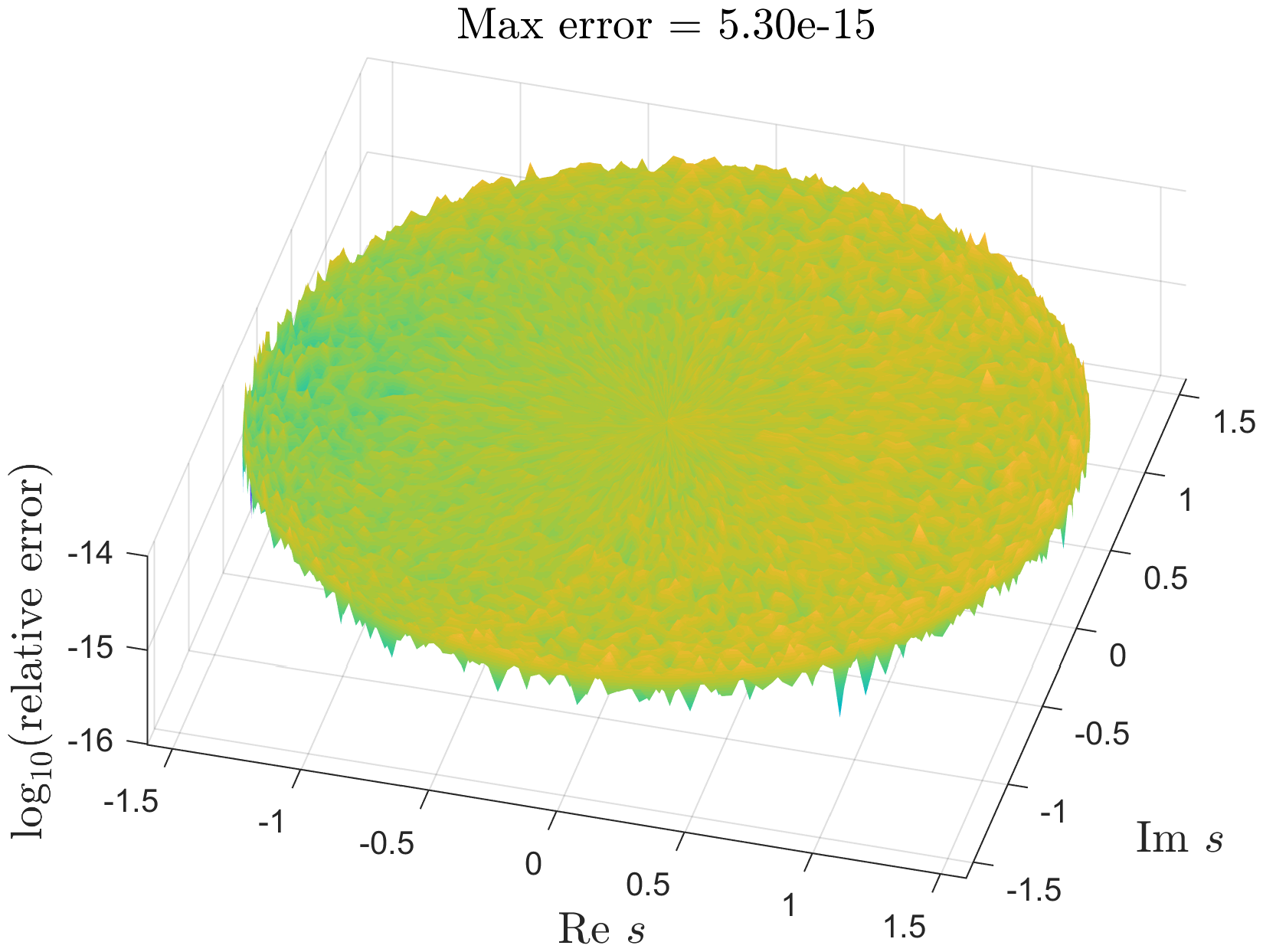}
 }
 \caption{
  The relative error on domain~I and domain~II in the $z$-plane (left) and on domain~III in the $s$-plane (right) for the approximation of $F(-1/3,1/2,1/2,z)$ computed with the multidomain spectral method. Recall that the $s$-plane corresponds to large $z$, the `far field', since $s = -1/(z-1/2)$.}
 \label{compF}
\end{figure}

The relative error is plotted in Figure~\ref{compF} for the test problem  in the $z$ and $s$ planes. Note that in the left frame that the error is largest close to the singular point $z = 1$ (due to the ill-conditioning of function evaluation in the vicinity of the singularity, as mentioned in the previous section). 

%

\section{Examples} \label{examples_sect}

In this section we consider further examples. The interesting paper 
\cite{POP} discussed challenging tasks for different numerical 
approaches and gave a table of 30 test cases for 5 different methods 
with recommendations when to use which. Note that our approach is 
complementary to \cite{POP}: we want to present an efficient approach 
to compute a solution to a Fuchsian equation, here the hypergeometric 
one, not for a single value, but on the whole compactified real line or on the 
whole Riemann sphere, and this for a wide range of the parameters 
$a,b,c$. To treat specific values of $a,b,c,z$, it is better to use the codes discussed in 
\cite{POP}. For generic values of the parameters, the present 
approach and the codes discussed in \cite{POP,PSPcode} produce similar 
results. 

The approach of this paper is supposed to fail when the 
genericness condition (\ref{generic}) is violated. Since we work 
with finite precision, the conditions (\ref{generic}) are supposed to 
hold for a whole range of parameters, i.e., that there are no 
integers in the intervals  $[c-\epsilon,c+\epsilon]$, 
$[c-a-b-\epsilon,c-a-b+\epsilon]$ and $[b-a-\epsilon,b-a+\epsilon]$ 
for some $\epsilon>0$. The used spectral methods are very sensitive 
to the possible appearance of logarithms in the solutions when condition
(\ref{generic}) does not hold, and thus there will be a loss of 
accuracy even in the vicinity of such cases. There will be either 
problems in the conditioning of the matrices (\ref{USsystem}) 
corresponding to the 5 ODEs introduced in section \ref{sec2}, or 
there will be problems with the matching conditions at the domain 
boundaries if no linearly independent solutions have been identified 
with necessary accuracy. 

Since most of the examples of \cite{POP} address degenerate or almost 
degenerate cases, they are outside the realm of applicability of the 
present approach. Below we present cases that can be treated with 
the present code together with additional examples along the lines of 
\cite{POP}. We define $\Delta F := 
|F_{num}(a,b,c,z)-F_{ex}(a,b,c,z)|/|F_{ex}(a,b,c,z)|$, where we use 
Maple with 30 digits as the  reference solution. 

We first 
address examples with real $z$ and give in Table~\ref{table1} the 
first 3 digits of the exact solutions, the quantity $\Delta F$  and 
the number of Chebychev coefficients $n$. It can be seen that a relative 
accuracy of the order of $10^{-10}$ can be reached even when the 
modulus of the hypergeometric function is of the order of $10^{-7}$. 

For the results in Table~\ref{table2} in which the argument $z$ is complex, the number of (i) Chebychev and (ii) Fourier coefficients of the solutions on the ellipses and (iii) the number Chebychev coefficients of the radial Fourier coefficients $u_k(r)$ are in the same ballpark as those required for the test problem (roughly $300$, $110$ and $40$ for (i), (ii) and (iii), respectively (see Figures~\ref{elI}, \ref{CFcoeffs} and \ref{elI2})).

\begin{table}[tbp]
    \centering
    {\renewcommand{\arraystretch}{1.25}
\renewcommand{\tabcolsep}{0.2cm}
   \begin{tabular}{|c|c|c|c|}
    \hline
    $a, b, c, z$ & $F(a,b,c,z)$ & $\Delta F$ & $n$  \\
    \hline
    $-0.1$, 0.2, 0.3, 0.5 & 0.956 & $1.2*10^{-16}$
     & 30  \\
     \hline
    $-0.1$, 0.2, 0.3, 1.5 &$0.904+0.179i$ & 
      $6.1*10^{-16}$ & 30  \\
   \hline
    $-0.1$, 0.2, 0.3, 100 & $1.365+0.400i$ & 
      $4.7*10^{-16}$ & 30  \\
    \hline
    $2 + 8i$, $3 -5i$,  $\sqrt{2} -i\pi$,  0.25 & $-3.670-4.764i$& 
    $7.9*10^{-15}$& 
    50\\
   \hline
    $2 + 8i$, $3 -5i$,  $\sqrt{2} -i\pi$,  0.75 &  
    $6882.463-6596.555i$& $8.3*10^{-15}$&
    50\\
     \hline
    $2 + 8i$, $3 -5i$,  $\sqrt{2} -i\pi$, $-10$ & $-0.0166-0.0067i$& 
    $7.5*10^{-15}$& 
    50\\
   \hline
    $2 + 200i$, $5 - 100i$, $10 + 500i$, 0.8 &$-4.103 + 6.013i$  & $5.9*10^{-15} $
    & 70 \\
    \hline
    2.25, 3.75, $-0.5$, $-1$ &  
    $-0.631$ & $4.3*10^{-12}$ & 50  \\
    \hline
     2 + 200i, 5, 10, 0.6 &  
     $(1.4997+5.771i)*10^{-7}$ &$2.4*10^{-10}$ & 160  \\
    \hline
   \end{tabular}}
   \caption{Examples  for the hypergeometric 
   function compared to a multiprecision computation in Maple for 
   real $z$.}
        \label{table1}
\end{table}

\begin{table}[tbp]
    \centering
    {\renewcommand{\arraystretch}{1.25}
\renewcommand{\tabcolsep}{0.2cm}
\begin{tabular}{|c|c|c|}
    \hline
    $a, b, c, z$ & $F(a,b,c,z)$ & $\Delta F$   \\
    \hline
    0.1, 0.2, $-0.3$, $-0.5 + 0.5i$ &  
    $1.027-0.013i$ & $2.3*10^{-16}$ \\
    \hline
    0.1, 0.2, $-0.3$, $1 + 0.5i$ &  
    $1.037-0.153i$ 
    & $6.4*10^{-16}$ \\
    \hline
    0.1, 0.2, $-0.3$, $5 + 5i$ &  
    $1.102+0.0288i$ 
    & $1.6*10^{-15}$ \\
    \hline
    4, 1.1, 2, $\exp(i\pi/3)$ & $ -0.461+0.487i$ & $4.0*10^{-14}$    \\
   \hline
    4, 1.1, 2, $1+5i$ & $-0.0183+0.0436i$ & 
    $9.1*10^{-14}$    \\
   \hline
    4, 1.1, 2, $-5+5i$ & $0.0216+0.0255i$ & 
    $9.1*10^{-14}$    \\
    \hline
    2/3, 1, 4/3, $\exp(i\pi/3)$ &  $0.883+0.50998i$ &$4.0*10^{-15}$   \\
   \hline
    2/3, 1, 4/3, $2i$ &  $0.562+0.373i$ &$7.1*10^{-15}$   \\
   \hline
    2/3, 1, 4/3, $1+i$ &  $0.740+0.740i$ &$4.5*10^{-15}$   \\
   \hline
    2/3, 1, 4/3, $100i$ &  $0.041+0.0609i$ &$8.7*10^{-15}$   \\
    \hline
  \end{tabular}}
   \caption{Examples  for the hypergeometric 
   function compared to a multiprecision computation in Maple for 
   complex $z$.}
        \label{table2}
\end{table}


\section{Outlook}
In this paper we have presented a spectral approach for the 
construction of the Gauss hypergeometric function on the whole 
Riemann sphere. One ingredient was essentially Kummer's approach to 
represent the solution to the hypergeometric function in the vicinity 
of each of the singularities 0, 1, $\infty$ via the hypergeometric 
function near 0. Since the transformation to obtain the second linearly 
independent solution to the hypergeometric equation near 0 is thus 
known, we did not address the task to compute also this solution.

The presented approach assumes a generic choice of the parameters 
$a,b,c$ for which no logarithms appear and for which the 
hypergeometric function is thus a function on a Riemann surface of 
finite genus. If the genericness condition (\ref{generic}) is not 
satisfied within a numerical precision of at least $10^{-6}$, the possible 
appearence of logarithms shows in the conditioning of the 
spectral matrices for the studied ODEs and the matching conditions where the 
matrix for the coefficients $\alpha,\beta$ (\ref{match1params})--(\ref{match1}) and for 
$\gamma,\delta$ (\ref{match2params})--(\ref{match2})
respectively can be singular. The latter implies that no linearly 
independent solutions have been identified. To address such 
cases, the following ansatz can be applied: if $y_{1}(x)$ is the
solution regular at $x=0$, then $y_{2}(x)=c_{0}y_{1}(x)\ln x + v(x)$ 
with $c_{0}=const$ can be a linearly 
independent solution to the hypergeometric equation, where 
$v(x)=x^{\kappa_{1}}\sum_{n=0}^{\infty}b_{n}x^{n}$ ($\kappa_{1}$ is 
one of the exponents of the symbol (\ref{symbol})). Thus $v$ 
satisfies an inhomogeneous Fuchsian equation which can be solved with 
a similar approach as before. To see whether such an ansatz allows for 
a similar accuracy for almost degenerate cases as for non-degenerate cases will be the subject of further 
research. 

One motivation of this work was to present an approach for general 
Fuchsian equations such as the Lam\'e and Heun equations. The latter equation represents a significant challenge with rich potential benefits---see for example \cite{FS12} for problems related to computation of the Heun function and its application to general relativity.  The 
main change here is the appearence of a fourth singularity which 
implies that a fourth domain needs to be introduced which in addition 
depends on a parameter. The rest of the 
approach remains unchanged. The  techniques used to study the
hypergeometric function as a meromorphic function on the Riemann 
sphere are also applicable to Painlev\'e transcendents as discussed 
in \cite{DubrovinGravaKlein,KS}. These nonlinear ordinary 
differential equations (ODEs) also have a wide range of applications, 
see \cite{Cla} and references therein. The similarity is due to the 
fact that Painlev\'e transcendents are meromorphic functions on the 
complex plane as is the case for the solutions of Fuchsian 
equations. Note that nonlinearities only affect the solution process 
on the real line and on the ellipses in the complex plane, i.e., 
one-dimensional problems. The only truely two-dimensional method, the 
solution of the Laplace equation for the interior of the ellipses, is 
unchanged for the Painlev\'e transcendents since the latter will be 
in general meromorphic as well. This replaces the task of solving a nonlinear ODE in the complex plane (which ultimately requires the solution of a system of nonlinear algebraic equations) with a linear PDE (which requires the solution of a linear system). The study of such transcendents, also 
on domains containing poles in the 
complex plane as in \cite{FW},  with the techniques 
outlined in this paper will be also subject to further research. 
Combining the compactification techniques of the present paper and 
the Pad\'e approach of \cite{FW}, it should be possible to study 
domains with a finite number of poles.

\section*{Acknowledgement}
This work was partially supported by the PARI and FEDER programs in 
2016 and 2017, by the ANR-FWF project ANuI and by the Marie-Curie 
RISE network IPaDEGAN. M.~Fasondini acknowledges financial support from the EPSRC grant EP/P026532/1.
 We thank C. Lubich for helpful remarks.

\end{document}